\documentclass[11pt]{amsart}
\usepackage{amscd,amssymb,latexsym,amsmath,cite,amsthm}
\usepackage{epsfig}
\def\eq#1{{\rm(\ref{#1})}}
\theoremstyle{plain}
\newtheorem{thm}{Theorem}[section]
\newtheorem{prop}[thm]{Proposition}
\newtheorem{lem}[thm]{Lemma}

\newtheorem{quest}[thm]{Question}
\theoremstyle{definition}
\newtheorem{dfn}[thm]{Definition}
\newtheorem{rem}[thm]{Remark}
\newtheorem{ans}[thm]{Ansatz}
\def\Re{\mathop{\rm Re}}
\def\Im{\mathop{\rm Im}}
\def\Image{\mathop{\rm Image}}
\def\SO{\mathbin{\rm SO}}
\def\U{\mathbin{\rm U}}
\def\vol{\mathop{\rm vol}}
\def\diag{\mathop{\rm diag}}
\def\ge{\geqslant}
\def\le{\leqslant}
\def\N{{\mathbin{\mathbb N}}}
\def\R{{\mathbin{\mathbb R}}}
\def\Z{{\mathbin{\mathbb Z}}}
\def\Q{{\mathbin{\mathbb Q}}}
\def\C{{\mathbin{\mathbb C}}}
\def\al{\alpha}
\def\be{\beta}
\def\ga{\gamma}

\def\io{\iota}
\def\ep{\epsilon}
\def\la{\lambda}
\def\th{\theta}
\def\De{\Delta}
\def\Si{\Sigma}
\def\pd{\partial}
\def\ts{\textstyle}

\def\iy{\infty}
\def\ra{\rightarrow}
\def\ab{\allowbreak}
\def\longra{\longrightarrow}
\def\t{\times}

\def\ti{\tilde}
\def\d{{\rm d}}
\def\sm{\setminus}
\def\ov{\overline}
\def\ha{{\ts\frac{1}{2}}}
\def\bs{\boldsymbol}
\def\md#1{\vert #1 \vert}
\def\bmd#1{\big\vert #1 \big\vert}
\def\ms#1{\vert #1 \vert^2}
\def\an#1{\langle #1 \rangle}

\begin{document}
\title[Self-similar solutions for Lagrangian MCF]{Self-similar
solutions and \\ translating solitons for \\  Lagrangian mean
curvature flow}
\author{Dominic Joyce, Yng-Ing Lee and Mao-Pei Tsui}
\address{The Mathematical Institute, 24-29 St. Giles,
Oxford, OX1 3LB, UK }
\email{joyce@maths.ox.ac.uk}
\address{Department of Mathematics and Taida Institute of
Mathematical Sciences, National Taiwan University, Taipei 10617,
Taiwan}
\email{yilee@math.ntu.edu.tw}
\address{Department of Mathematics, University of Toledo,
Toledo, OH 43606, USA, and  Taida Institute of Mathematical
Sciences, Taipei, Taiwan}
\email{email: mao-pei.tsui@utoledo.edu}
\begin{abstract}
We construct many self-similar and translating solitons for
Lagrangian mean curvature flow, including self-expanders and
translating solitons with arbitrarily small oscillation on the
Lagrangian angle. Our translating solitons play the same role as
cigar solitons in Ricci flow, and are important in studying the
regularity of Lagrangian mean curvature flow.

Given two transverse Lagrangian planes $\R^n$ in $\C^n$ with sum of
characteristic angles less than $\pi$, we show there exists a
Lagrangian self-expander asymptotic to this pair of planes. The
Maslov class of these self-expanders is zero. Thus they can serve as
local models for surgeries on Lagrangian mean curvature flow.
Families of self-shrinkers and self-expanders with different
topologies are also constructed. This paper generalizes the work of
Anciaux \cite{Anci}, Joyce \cite{Joyc1}, Lawlor \cite{Lawl} and Lee
and Wang~\cite{LeWa1,LeWa2}.
\end{abstract}
\maketitle

\section{Introduction}
\label{lm1}

Special Lagrangian submanifolds in Calabi--Yau $n$-folds have
received much attention in recent years, as they are key ingredients
in the Strom\-in\-ger--Yau--Zaslow Conjecture \cite{SYZ}, which
explains Mirror Symmetry of Calabi-Yau 3-folds. Thomas and Yau
\cite{ThYa} defined a notion of stability for graded Lagrangians $L$
in a Calabi--Yau $n$-fold $M$, and conjectured that if $L$ is stable
then the {\it Lagrangian mean curvature flow} of $L$ exists for all
time and converges to a special Lagrangian submanifold $L_\iy$ in
$M$, which should be the unique special Lagrangian in the
Hamiltonian equivalence class of~$L$.

Rewriting this in terms of the derived Fukaya category $D^b{\rm
Fuk}(M)$ of $M$, as in Kontsevich's Homological Mirror Symmetry
programme \cite{Kont}, and using Bridgeland's notion of stability
condition on triangulated categories \cite{Brid}, one can state an
improved (but still over-simplified) version of the Thomas--Yau
conjecture as follows: for any Calabi--Yau $n$-fold $M$, there
should exist a Bridgeland stability condition $(Z,{\mathcal P})$ on
$D^b{\rm Fuk}(M)$ depending on the holomorphic $(n,0)$-form $\Omega$
on $M$, such that a graded Lagrangian $L$ in $M$ is $(Z,{\mathcal
P})$-stable, regarded as an object in $D^b{\rm Fuk}(M)$, if and only
if the Lagrangian mean curvature flow of $L$ exists for all time and
converges to a special Lagrangian submanifold $L_\iy$ in $M$, which
should be unique in the isomorphism class of $L$ in $D^b{\rm
Fuk}(M)$. A related method for constructing special Lagrangians by
minimizing volume amongst Lagrangians using Geometric Measure Theory
was proposed by Schoen and Wolfson~\cite{ScWo}.

To carry these programmes through to their conclusion will require
a deep understanding of Lagrangian mean curvature flow, and of the
possible singularities that can occur during it in finite time.
Singularities in Lagrangian mean curvature flow are generally
locally modelled on {\it soliton solutions}, such as Lagrangians
in $\C^n$ which are moved by rescaling or translation by mean
curvature flow. There are two important results in this area. The
first one is due to Wang \cite{Wang},  who observed that mean
curvature flow for almost calibrated Lagrangians in Calabi--Yau
$n$-folds cannot develop type I singularities.   And the second
one is due to Neves \cite{Neve}, who (loosely) proved that
singularities of such flows are modelled to leading order on
special Lagrangian cones when applying central blow up near the
singularities.

In this paper, we construct many examples of self-similar solutions
and translating solitons for Lagrangian mean curvature flow. Our
Lagrangians $L$ in $\C^n$ are the total space of a 1-parameter
family $Q_s$, $s\in I$, where $I$ is an open interval in $\R$, and
each $Q_s$ is a quadric in a Lagrangian plane $\R^n$ in $\C^n$,
which evolve according to an o.d.e.\ in $s$. The construction
includes and generalizes examples of Lagrangian solitons or special
Lagrangians due to Anciaux \cite{Anci}, Joyce \cite{Joyc1}, Lawlor
\cite{Lawl}, and Lee and Wang~\cite{LeWa1,LeWa2}.

The authors believe that two of our families of examples may have
particular significance for future work on Lagrangian mean curvature
flow. Firstly, in Theorems C and D of \S\ref{lm32}, we show that if
$L_1,L_2$ are transverse Lagrangian planes in $\C^n$ and the sum of
characteristic angles of $L_1,L_2$ is less than $\pi$, and $\al>0$,
then we can construct a unique closed, embedded Lagrangian
self-expander $L$ with rate $\al$ diffeomorphic to ${\mathcal
S}^{n-1}\t\R$ and asymptotic to $L_1\cup L_2$. These examples could
be used as local models for surgeries during Lagrangian mean
curvature flow.

As in the Ricci flow proof of the Poincar\'e conjecture \cite{Pere},
it seems likely that to get long-time existence for Lagrangian mean
curvature flow, it will be necessary to allow the flow to develop
singularities, and continue the flow after a surgery which changes
the topology of the Lagrangian. Research by the first author
(unpublished) indicates that an important condition in the improved
Thomas--Yau conjecture described above is that the Lagrangians
should have {\it unobstructed Lagrangian Floer homology}, in the
sense of Fukaya, Oh, Ohta and Ono \cite{FOOO, FOO1}. But mean
curvature flow amongst nonsingular, immersed Lagrangians can cross
`walls', on the other side of which Lagrangian Floer homology is
obstructed. When this happens, the correct thing to do is to do a
surgery, and glue in a Lagrangian self-expander from Theorems C
and~D.

Secondly, in Corollary I of \S\ref{lm34} we give an explicit family
of closed, embedded Lagrangian translating solitons $L$ in $\C^n$
for $n\ge 2$, which are diffeomorphic to $\R^n$, and asymptotic in a
weak sense to a union $L_1\cup L_2$ of Lagrangian planes
$L_1,L_2\cong\R^n$ in $\C^n$, with $L_1\cap L_2\cong\R$. The
oscillation of the Lagrangian angle of $L$ can be chosen arbitrarily
small. If these examples can arise as local models for finite time
singularities for Lagrangian mean curvature flow, they may represent
a kind of bad behaviour, which could cause difficulties with the
Thomas--Yau programme even in dimension~2.

As well as these two families, we construct new examples of compact,
immersed Lagrangian self-shrinkers in $\C^n$ diffeomorphic to ${\mathcal
S}^1\t{\mathcal S}^{n-1}$, of closed, immersed Lagrangian self-expanders
and self-shrinkers diffeomorphic to ${\mathcal S}^1\t{\mathcal S}^{m-1}\t
\R^{n-m}$ for $0<m<n$, of non-closed, immersed Lagrangian
self-expanders diffeomorphic to ${\mathcal S}^{m}\t\R^{n-m}$ for
$0<m<n-1$, of non-closed, immersed Lagrangian self-shrinkers
diffeomorphic to ${\mathcal S}^{m}\t\R^{n-m}$ for $0<m<n$, and of
closed, embedded Lagrangian translating solitons diffeomorphic to
$\R^n$ with infinite oscillation of the Lagrangian angle. These
examples include those of Anciaux \cite{Anci} and Lee and Wang
\cite{LeWa1,LeWa2}, and approach the special Lagrangians of Joyce
\cite{Joyc1} and Lawlor \cite{Lawl} in a limit.

We begin in \S\ref{lm2} with some background material. Our main
results are stated and discussed in \S\ref{lm3}, which is the part
of the paper we intend most people to actually read. The proofs are
given in \S\ref{lm4}--\S\ref{lm6}.
\smallskip

\noindent{\it Acknowledgements.} Part of this paper was completed
while the authors were visiting Taida Institute of Mathematical
Sciences  (TIMS) in National Taiwan University, Taipei, Taiwan. The
authors wish to express their gratitude for the excellent support
provided by the centre during their stays. The first author would
like to thank Richard Thomas and Tom Ilmanen for useful
conversations. The second author would like to express her special
gratitude to  Mu-Tao Wang for many enlightening discussions and
collaborations on the subject over the years, which benefitted her a
lot. The calculation techniques used in the proof of Theorem A were
first observed in \cite{LeWa1,LeWa2}. She would also like to thank
Andr\'e Neves for his interests in this work and enlightening
discussions. The third author would like to thank Duong Hong Phong,
Shing-Tung Yau and Mu-Tao Wang for their constant advice,
encouragement and support.

\section{Background material}
\label{lm2}

Our ambient space is always the complex Euclidean space $\C^n$ with
coordinates $z_j=x_j+iy_j$, the standard symplectic form
$\omega=\sum_{j=1}^n\d x_j\wedge \d y_j$, and the standard almost
complex structure $J$ with $J(\frac{\partial}{\partial
x_j})=\frac{\partial}{\partial y_j}$. A {\it Lagrangian
submanifold\/} is an $n$-dimensional submanifold in $\C^n$ on which
the symplectic form $\omega$ vanishes. On a Lagrangian submanifold
$L$, the mean curvature vector ${H}$ is given by
\begin{equation}
{H}=J\nabla \th,
\label{lm2eq1}
\end{equation}
where $\th$ is the {\it Lagrangian angle\/} and $\nabla$ is the
gradient on $L$. The angle function $\th:L\ra\R$ or
$\th:L\ra\R/2\pi\Z$ can be defined by the relation that $\d
z_1\wedge\cdots\wedge \d z_n\vert_L\equiv e^{i\th}\vol_L$. When
$\cos\theta\geq\epsilon$ on $L$ for some positive $ \epsilon > 0$,
$L$ is called {\it almost-calibrated}. The {\it Maslov class} on $L$
is defined by the cohomology class of $\d\theta$. Hence $L$ is
Maslov zero when $\th$ is a globally defined function from $L$ to
$\R$.

By the first variation formula, the {\it mean curvature vector\/}
points in the direction in which the volume decreases most rapidly.
{\it Mean curvature flow\/} deforms the submanifold in the direction
of the mean curvature vector. As special Lagrangians are volume
minimizing, it is natural to use mean curvature flow to construct
special Lagrangians. Equation \eq{lm2eq1} implies that mean
curvature flow is a Lagrangian deformation, that is, a Lagrangian
submanifold remains Lagrangian under mean curvature flow, as in
Smoczyk~\cite{Smoc}.

A Lagrangian submanifold $L$ in $\C^n$ is fixed by mean curvature
flow if and only if the Lagrangian angle $\th$ on $L$ is constant,
that is, if and only if $L$ is {\it special Lagrangian\/} with
phase $e^{i\th}$, as in Harvey and Lawson \cite[\S III]{HaLa}. A
Lagrangian $L$ in $\C^n$ is called {\it Hamiltonian stationary\/}
if the Lagrangian angle $\th$ on $L$ is harmonic, that is, if
$\De\th=0$ on $L$. This implies that the volume of $L$ is
stationary under Hamiltonian  deformations.

In geometric flows such as Ricci flow or mean curvature flow,
singularities are often locally modelled on soliton solutions. In
the case of mean curvature flows, two types of soliton solutions
of particular interest are those moved by scaling or translation
in Euclidean space. We recall that solitons moved by scaling must
be of the form:

\begin{dfn} A submanifold $L$ in Euclidean space $\R^n$ is
called a {\it self-similar solution\/} if $H\equiv\al F^\perp$ on
$L$ for some constant $\al$ in $\R$, where $F^\perp$ is the
projection of the position vector $F$ in $\R^n$ to the normal bundle
of $L$, and $H$ is the mean curvature vector of $L$ in $\R^n$. It is
called a {\it self-shrinker\/} if $\al<0$ and {\it self-expander\/}
if~$\al>0$.
\label{lm2def1}
\end{dfn}

It is not hard to see that if $F$ is a self-similar solution, then
$F_t$ defined by $F_t=\sqrt{2\al t}\,F$ is moved by the mean
curvature flow. By Huisken's monotonicity formula \cite{Huis}, any
central blow up of a finite-time singularity of the mean curvature
flow is a self-similar solution. When $\al=0$, the submanifold is
minimal. The submanifolds which are moved by translation along
mean curvature flow must be of the form:

\begin{dfn} A submanifold $L$ in Euclidean space $\R^n$ is called a
{\it translating soliton\/} if there exists a constant vector $T$ in
$\R^n$ such that $H+V\equiv T$ on $L$, where $V$ is the component of
$T$ tangent to $L$, and $H$ is the mean curvature vector of $L$ in
$\R^n$. An equivalent equation is $H\equiv T^\perp$. The 1-parameter
family of submanifolds $L_t$ defined by $L_t=L+t\,T$ for $t\in\R$ is
then a solution to mean curvature flow, and we call $T$ a {\it
translating vector}.
\label{lm2def2}
\end{dfn}

\begin{dfn} A translating soliton is called a {\it gradient
translating soliton} if $V=\nabla f$ for some smooth function
$f:L\ra\R$.
\label{lm2def3}
\end{dfn}

Any translating soliton for mean curvature flow in $\R^n$ must be
a gradient translating soliton. Since this simple fact does not
appear in the literature, we include a proof here for
completeness.

\begin{prop} A translating soliton in $\R^n$ that satisfies $H+V=T$
where $T$ is a constant vector must be a gradient translating
soliton. In fact, $H+\nabla\an{T,F}=T,$ where $F$ is the position
vector.
\label{lm2prop1}
\end{prop}

\begin{proof} Let $F=(F_1,\ldots,F_n)$ be the position function
and $\ov{\nabla}$ be the standard connection in $\R^n$. Then
$\ov{\nabla}F_i=e_i$. We may write $T=\sum_{i=1}^nT^ie_i=
\ov{\nabla}\an{T, F}$. Then $(\ov{\nabla}\an{T,F})^\top=V $ and
$\nabla\an{T,F}\vert_L=V$, where $(\ov{\nabla}\an{T,F})^\top$ is the
orthogonal projection of $\ov{\nabla}\an{T,F}$ to $TL$. This shows
$L$ is a gradient translating soliton.
\end{proof}

Here is a counterpart of this result for {\it Lagrangian}
translating solitons.

\begin{prop} A connected Lagrangian $L$  in
$\C^n$ is a translating soliton with translating vector $T$ if and
only if\/ $\th\equiv -\an{JT,F}\vert_L+c$ for some $c\in\R,$ where
$F$ is the position vector. Thus a Lagrangian translating soliton is
Maslov zero. \label{lm2prop2}
\end{prop}

\begin{proof} Suppose $L$ is a translating soliton with translating
vector $T$. We have $J\,\nabla\th\equiv H\equiv T^\perp$ as
sections of the normal bundle of $L$ in $\C^n$. Applying $-J$
gives $\nabla\th\equiv -J\bigl(T^\perp)=-(JT)^\top$, as sections
of the tangent bundle of $L$. We then have
$\nabla\th=-\nabla\an{JT,F}$ as in the proof of
Proposition~\ref{lm2prop1}.  Because $L$ is connected, it follows
that  $\th\equiv -\an{JT,F}\vert_L+c$ for some $c\in\R$. Hence
$\th$ can be lifted from $\R/2\pi\Z$ to $\R$, and $L$ is Maslov
zero. Conversely, suppose $\th\equiv -\an{JT,F}\vert_L+c$. Then
$\nabla\th=-\nabla\an{JT,F}=-(\ov{\nabla}\an{JT,F})^\top=
-(JT)^\top.$ It follows that $H=J\nabla\th=T^\perp$, and $L$ is a
translating soliton with vector~$T$.
\end{proof}

\section{Statements of main results}
\label{lm3}

We now state and briefly discuss the main results of this paper. The
results of \S\ref{lm31}, \S\ref{lm32} and \S\ref{lm33} will be
proved in sections \ref{lm4}, \ref{lm5} and \ref{lm6}, respectively.
The proofs of results in \S\ref{lm34} are brief, and are included
there.

\subsection{An ansatz for self-similar Lagrangians}
\label{lm31}

The following ansatz describes the class of $n$-submanifolds of
$\C^n$ amongst which we will seek examples self-similar solutions
for Lagrangian mean curvature flow.

\begin{ans} Fix $n\ge 1$. Consider $n$-submanifolds $L$ in $\C^n$ of
the form:
\begin{equation}
L=\bigl\{\bigl(x_1w_1(s),\ldots,x_nw_n(s)\bigr):\text{$s\in I$,
$x_j\in\R$, $\ts\sum_{j=1}^{n}\la_jx_j^2=C$}\bigl\},
\label{lm3eq1}
\end{equation}
where $\la_1,\ldots,\la_n,C\in\R\sm\{0\}$ are nonzero constants, $I$
is an open interval in $\R$, and $w_1,\ldots,w_n:I\ra\C\sm\{0\}$ are
smooth. We want $L$ to satisfy:
\begin{itemize}
\setlength{\itemsep}{0pt}
\setlength{\parsep}{0pt}
\item[(i)] $L$ is Lagrangian;
\item[(ii)] the Lagrangian angle $\th:L\ra\R$ or
$\th:L\ra\R/2\pi\Z$ of $L$ is a function only of $s$, not of
$x_1,\ldots,x_n$; and
\item[(iii)] $L$ is a self-similar solution under mean curvature
flow in $\C^n$, that is, $H\equiv\al F^\perp$ on $L$ as in
Definition~\ref{lm2def1}.
\end{itemize}
\label{lm3ans1}
\end{ans}

The motivation for this ansatz is that it includes, and generalizes,
several families of examples in the literature. For special
Lagrangian submanifolds, with $\th\equiv 0$ in (ii) and $\al=0$ in
(iii), the ansatz includes the examples of Lawlor \cite{Lawl} with
$\la_1=\cdots=\la_n=C=1$, and Joyce \cite[\S 5--\S 6]{Joyc1}. For
more general Lagrangian self-similar solutions, it includes the
examples of Abresch and Langer \cite{AbLa} when $n=1$, the examples
of Anciaux \cite{Anci}, which have $\la_1=\cdots=\la_n=C=1$ and are
symmetric under the action of $\SO(n)$ on $\C^n$, and the examples
of Lee and Wang \cite[\S 6]{LeWa1}, \cite{LeWa2}, which have
$w_j(s)\equiv e^{i\la_js}$. M.-T. Wang and the second author also
tried to study an ansatz of a similar form to \eq{lm3eq1} before.

It is a long but straightforward calculation to find the conditions
on $\la_1,\ab\ldots,\ab\la_n,\ab C,w_1,\ldots,w_n$ for $L$ in
\eq{lm3eq1} to be Lagrangian, to compute its Lagrangian angle $\th$
and mean curvature $H$, and to work out whether $L$ is a
self-similar solution to Lagrangian mean curvature flow. In this way
we prove the following theorem.
\medskip

{\bf Theorem A.} {\it Let\/ $\la_1,\ldots,\la_n,C\in\R\sm\{0\}$
and\/ $\al\in\R$ be constants, $I$ be an open interval in $\R,$
and\/ $\th:I\ra\R$ or $\th:I\ra\R/2\pi\Z$ and\/
$w_1,\ldots,w_n:I\ra\C\sm\{0\}$ be smooth functions. Suppose that
\begin{equation}
\begin{aligned}
\frac{\d w_j}{\d s}&= \la_je^{i\th(s)}\,\ov{w_1\cdots w_{j-1}
w_{j+1}\cdots w_n},\qquad j=1,\ldots,n,\\
\frac{\d\th}{\d s}&=\al\Im(e^{-i \th(s)} w_1\cdots w_n),
\end{aligned}
\label{lm3eq2}
\end{equation}
hold in $I$. Then the submanifold\/ $L$ in $\C^n$ given by
\begin{equation}
L=\bigl\{\bigl(x_1w_1(s),\ldots,x_nw_n(s)\bigr):\text{$s\in I,$
$x_j\in\R,$ $\ts\sum_{j=1}^{n}\la_jx_j^2=C$}\bigl\},
\label{lm3eq3}
\end{equation}
is Lagrangian, with Lagrangian angle $\th(s)$ at\/
$(x_1w_1(s),\ldots,x_nw_n(s)),$ and its position vector\/ $F$ and
mean curvature vector\/ $H$ satisfy\/ $\al F^\perp=CH$. That is,
$L$ is a self-expander when\/ $\al/C>0$ and a self-shrinker when\/
$\al/C<0$. When $\al=0$ the Lagrangian angle $\th$ is constant, so
that\/ $L$ is special Lagrangian, with\/ $H=0$. In this case the
construction reduces to that of Joyce~{\rm\cite[\S 5]{Joyc1}}.}
\medskip

We can simplify the equations \eq{lm3eq2}, generalizing~\cite[\S
5.2]{Joyc1}.
\medskip

{\bf Theorem B.} {\it In the situation of Theorem A, let\/
$w_1,\ldots,w_n,\th$ satisfy\/ \eq{lm3eq2}. Write $w_j\equiv
r_je^{i\phi_j}$ and\/ $\phi=\sum_{j=1}^n\phi_j$, for functions
$r_j:I\ra(0,\iy)$ and\/ $\phi_1,\ldots,\phi_n,\phi:I\ra\R$ or
$\R/2\pi\Z$. Fix $s_0\in I$. Define $u:I\ra\R$ by
\begin{equation*}
u(s)=2\int_{s_0}^sr_1(t)\cdots
r_n(t)\cos\bigl(\phi(t)-\th(t)\bigr)\d t.
\end{equation*}
Then $r_j^2(s)\equiv\al_j+\la_ju(s)$ for $j=1,\ldots,n$ and\/ $s\in
I,$ where $\al_j=r_j^2(s_0)$. Define a degree $n$ polynomial $Q(u)$
by $Q(u)=\prod_{j=1}^n(\al_j+\la_ju)$. Then the system of equations
\eq{lm3eq2} can be rewritten as
\begin{equation}
\begin{cases}
\begin{aligned}
\frac{\d u}{\d s}&=2Q(u)^{1/2}\cos(\phi-\th),\\
\frac{\d\phi_j}{\d s} &= -\frac{\la_jQ(u)^{1/2}\sin(\phi-\th)}{
\al_j+\la_ju}\,,\qquad j=1,\ldots,n,\\
\frac{\d\phi}{\d s}&=-Q(u)^{1/2}(\ln Q(u))' \sin(\phi-\th),\\
\frac{\d\th}{\d s} &= \al Q(u)^{1/2}\sin(\phi-\th).
\end{aligned}
\end{cases}
\label{lm3eq4}
\end{equation}
The Lagrangian self-similar solution $L$ in Theorem A may be
rewritten
\begin{equation}
\begin{split}
L=\bigl\{(x_1\sqrt{\al_1+\la_1u(s)}\,&e^{i\phi_1(s)},
\ldots,x_n\sqrt{\al_n+\la_nu(s)}\,e^{i\phi_n(s)}):\\
&\text{$x_1,\ldots,x_n\in\R,$ $s\in I,$
$\ts\sum_{j=1}^n\la_jx_j^2=C$}\bigr\}.
\end{split}
\label{lm3eq5}
\end{equation}
Moreover, for some $A\in\R$ the equations \eq{lm3eq4} have the first
integral}
\begin{equation}
Q(u)^{1/2}e^{\al u/2}\sin(\phi-\th)\equiv A.
\label{lm3eq6}
\end{equation}

\begin{rem} There is a lot of freedom to rescale the constants in
Theorems A and B without changing the Lagrangian $L$. In particular:\\
{\bf (a)} Set
\begin{gather*}
\ti\la_j=C\la_j/\md{C\la_j},\;\> \ti C=1,\;\> \ti\al=\al/C,\;\>
\ti I=C\md{C}^{-n/2}\ts\prod_{j=1}^n\md{\la_j}^{1/2}\cdot I,\\
\ti s=C\md{C}^{-n/2}\ts\prod_{j=1}^n\md{\la_j}^{1/2}s,\;\>
\ti\th=\th,\;\> \ti w_j=\md{C}^{1/2}\md{\la_j}^{-1/2}w_j,\\
\ti x_j=\md{C}^{-1/2}\md{\la_j}^{1/2}x_j,\;\> \ti
r_j=\md{C}^{1/2}\md{\la_j}^{-1/2}r_j,\;\>
\ti\phi_j=\phi_j,\;\> \ti\phi=\phi,\\
\ti u=Cu,\;\> \ti\al_j=\md{C}\md{\la_j}^{-1}\al_j,\;\> \ti
A=\md{C}^{n/2}\ts\prod_{j=1}^n\md{\la_j}^{-1/2}A,
\end{gather*}
where we regard $\ti w_j,\ti\th,\ti r_j,\ti\phi_j,\ti\phi,\ti u$ as
functions of $\ti s$ rather than $s$, so that $\ti w_j(\ti s)=
\md{C}^{1/2}\md{\la_j}^{-1/2}w_j(C^{-1}\md{C}^{n/2}
\prod_{j=1}^n\md{\la_j}^{-1/2}\ti s)$, for instance. Then these
$\ti\la_1,\ldots,\ti A$ satisfy Theorems A and B, with the same
Lagrangian $L$. Thus without loss of generality we can suppose
$\la_j=\pm 1$ for all $j$ and~$C=1$.

\noindent{\bf(b)} Translations of $I$ in $\R$, so that $I\mapsto
I+c$, $s\mapsto s+c$ also do not change $L$. Thus we can fix $0\in
I$ and $s_0=0$ in Theorem B.

\noindent{\bf(c)} Changing $s\mapsto -s$ and $\th\mapsto\th+\pi$
gives a solution with the same $L$, but the opposite orientation.

\noindent{\bf(d)} Changing $w_1\mapsto -w_1$, $s\mapsto -s$,
$\phi_1\mapsto \phi_1+\pi$, $\phi\mapsto\phi+\pi$, and $A\mapsto
-A$ gives a solution with the same $L$.  When $A=0$, as
$Q(u)=r_1\cdots r_n>0$ equation \eq{lm3eq6} gives
$\sin(\phi-\th)\equiv 0$, so \eq{lm3eq4} implies that
$\phi_j,\phi,\th$ are constant. Thus $L$ is an open subset of the
special Lagrangian $n$-plane
\begin{equation*}
\bigl\{(e^{i\phi_1}x_1,\ldots,e^{i\phi_n}x_n):x_1,\ldots,x_n\in\R\bigr\}
\end{equation*}
in $\C^n$. As we are not interested in this case, we will either
take~$A>0$ or $A<0$. For the case of explicit Lagrangian
self-expanders in \S\ref{lm5}, we choose $A<0$ to make the solutions
have a similar expression to Lawlor's examples, while we take $A>0$
in~\S\ref{lm6}.

\noindent{\bf(e)} If $t>0$, changing $\al\mapsto t^{-2}\al$,
$I\mapsto t^{n-2}I$, $s\mapsto t^{n-2}s$, $w_j\mapsto tw_j$,
$r_j\mapsto tr_j$, $u\mapsto t^2u$, $\al_j\mapsto t^2\al_j$,
$A\mapsto t^{n}A$ gives another solution for the Lagrangian $tL$
rather than $L$. Thus, by allowing rescalings $L\mapsto tL$, we
can also set $\al$ to $1,-1$ or 0. But we shall retain the
parameter $\al$, since taking limits $\al\ra 0$ shows how our
Lagrangian self-expanders or self-shrinkers are related to special
Lagrangian examples. \label{lm3rem1}
\end{rem}

\subsection{A class of explicit Lagrangian self-expanders}
\label{lm32}

As in Remark \ref{lm3rem1}, in Theorems A and B we may without
loss of generality suppose that $\la_j=\pm 1$, $C=1$ and $A<0$. We
now consider the case in which $\la_1=\cdots=\la_n=C=1$, $\al\ge
0$ and $A<0$. Then the Lagrangians $L$ we get are embedded and
diffeomorphic to ${\mathcal S}^{n-1}\t\R$. When $\al=0$ they are
the special Lagrangian `Lawlor necks' found by Lawlor \cite{Lawl}
and studied by Harvey \cite[p.~139--143]{Harv}, and Theorem C
below generalizes Harvey's treatment. For $\al>0$ they are
Lagrangian self-expanders. When $\al>0$, $a_1=\cdots=a_n$ and
$\psi_1=\cdots=\psi_n$ in Theorem C, the self-expander $L$ is
invariant under $\SO(n)$, and is one of the examples found by
Anciaux~\cite{Anci}.
\medskip

{\bf Theorem C.} {\it In Theorems A and B, suppose that\/
$\la_1=\cdots=\la_n=C=1$, $\al\ge 0$ and\/ $A<0$. Then any
solution of\/ {\rm\eq{lm3eq2},} or equivalently of\/
{\rm\eq{lm3eq4},} on an interval\/ $I$ in $\R$ can be extended to
a unique largest open interval\/ $I_{\max}$ in $\R$. Take
$I=I_{\max}$. Then by changing variables from $s$ in $I_{\max}$ to
$y=y(s)$ in $\R,$ we may rewrite the Lagrangian self-expander $L$
of\/ \eq{lm3eq3} and\/ \eq{lm3eq5} explicitly as follows.
Conversely, every $L$ of the following form comes from Theorems A
and B with\/ $\la_1=\cdots=\la_n=C=1$, $\al\ge 0$ and\/~$A<0$.

For given constants $\al\ge 0,$ $a_1,\ldots,a_n>0$ and\/
$\psi_1,\ldots,\psi_n\in\R,$ define $w_j(y)=e^{i\phi_j(y)}r_j(y)$
for $j=1,\ldots,n$ and\/ $y\in\R$ by
\begin{equation*}
r_j(y)=\sqrt{\ts\frac{1}{a_j}+y^2} \quad\text{and}\quad
\phi_j(y)=\psi_j+\int_0^y \frac{\d
t}{(\frac{1}{a_j}+t^2)\sqrt{P(t)}}\,,
\end{equation*}
where\/ $P(t)=\ts\frac{1}{t^2}\bigl(\prod_{k=1}^n(1+a_kt^2)e^{\al
t^2}-1\bigl)$. Then
\begin{equation}
L=\bigl\{\bigl(x_1w_1(y),\ldots,x_n w_n(y)\bigr):
\text{$x_1,\ldots,x_n\in\R,$ $\ts\sum_{j=1}^nx_j^2=1$}\bigr\}
\label{lm3eq7}
\end{equation}
is a closed, embedded Lagrangian diffeomorphic to ${\mathcal
S}^{n-1}\t\R$ and satisfying\/ $\al F^\perp=H$. If\/ $\al>0$ it is a
self-expander, and if\/ $\al=0$ it is one of Lawlor's examples of
special Lagrangian submanifolds {\rm\cite{Lawl}}. It has Lagrangian
angle}
\begin{equation}
\th(y)=\ts\sum_{j=1}^n\phi_j(y)+\arg\bigl(y+iP(y)^{-1/2}\bigr).
\label{lm3eq8}
\end{equation}

We can describe the asymptotic behaviour of these Lagrangians:
\medskip

{\bf Theorem D.} {\it In the situation of Theorem C, there exist\/
$\bar\phi_1,\ldots,\bar\phi_n\in(0,\frac{\pi}{2}]$ with\/
$\bar\phi_j=\int_0^\infty \frac{\d
t}{(\frac{1}{a_j}+t^2)\sqrt{P(t)}}$ for\/ $ j=1,\ldots,n,$ such that
the Lagrangian $L$ is asymptotic at infinity to the union of
Lagrangian planes $L_1 \cup L_2,$ where
\begin{align*}
L_1&=\bigl\{(e^{i(\psi_1+\bar{\phi}_1)}t_1,\ldots,
e^{i(\psi_n+\bar{\phi}_n)}t_n):t_1,\ldots,t_n\in\R\bigr\},\\
L_2&=\bigl\{(e^{i(\psi_1-\bar{\phi}_1)}t_1,\ldots,
e^{i(\psi_n-\bar{\phi}_n)}t_n):t_1,\ldots,t_n\in\R\bigr\}.
\end{align*}
We have $0<\bar\phi_1+\cdots+\bar\phi_n<\frac{\pi}{2}$ if\/ $\al>0,$
and\/ $\bar\phi_1+\cdots+\bar\phi_n=\frac{\pi}{2}$ if\/~$\al=0$.

Fix $\al>0$. Then $\Phi^n:(a_1,\ldots,a_n)\mapsto(\bar\phi_1,
\ldots,\bar\phi_n)$ gives a diffeomorphism
\begin{equation*}
\Phi^n:(0,\iy)^n\longra\bigl\{(\bar\phi_1,\ldots,\bar\phi_n)\in(0,
\ts\frac{\pi}{2})^n:
0<\bar\phi_1+\cdots+\bar\phi_n<\frac{\pi}{2}\bigr\}.
\end{equation*}
That is, for all\/ $\al>0$ and\/ $L_1,L_2$ satisfying
$0<\bar\phi_1+\cdots+ \bar\phi_n<\frac{\pi}{2}$ as above, Theorem C
gives a unique Lagrangian expander $L$ asymptotic to~$L_1\cup L_2$.

When $\al=0,$ it is studied by  Lawlor in {\rm \cite{Lawl}}. The map
$\Phi^n:(a_1,\ldots,a_n)\allowbreak\mapsto(\bar\phi_1,
\ldots,\bar\phi_n)$ gives a surjection
\begin{equation*}
\Phi^n:(0,\iy)^n\longra\bigl\{(\bar\phi_1,\ldots,\bar\phi_n)\in(0,
\ts\frac{\pi}{2}]^n:
\bar\phi_1+\cdots+\bar\phi_n=\frac{\pi}{2}\bigr\},
\end{equation*}
such that\/ $(a_1,\ldots,a_n)$ and\/ $(a_1',\ldots,a_n')$ have the
same image $(\bar\phi_1,\ldots,\bar\phi_n)$ if and only if\/
$a_j'=ta_j$ for some $t>0$ and all\/ $j=1,\ldots,n,$ and the
corresponding special Lagrangians $L,L'$ satisfy~$L'=t^{-1/2}L$. }

\medskip

By applying an element of $\U(n)$, Theorem D also shows that we
can construct a unique Lagrangian self-expander with constant
$\al$ asymptotic to any pair of Lagrangian planes in $\C^n$ which
intersect transversely at the origin and have sum of
characteristic angles less than $\pi$. As the union of a pair of
planes is volume minimizing if and only if the sum of
characteristic angles is greater or equal to $\pi$ \cite{Lawl},
our result is sharp.

The Lagrangian self-expanders in Theorems C and D have {\it
arbitrarily small oscillation of the Lagrangian angle}. That is, if
$\sum_{j=1}^n \bar{\phi}_j=\frac{\pi}{2}-\ep$ in Theorem D, then
\eq{lm3eq8} implies that $L$ in \eq{lm3eq7} has Lagrangian angle
varying in $(\sum_{j=1}^n\psi_j+\frac{\pi}{2}-\ep,\sum_{j=1}^n\psi_j
+\frac{\pi}{2}+\ep)$, an open interval of width $2\ep$, which can be
made arbitrarily small. Thus these self-expanders are {\it almost
calibrated}, and Maslov zero.

The `Lawlor necks' \cite{Lawl} have been used as local models in
resolving intersection points of special Lagrangians, see for
example Butscher \cite{Buts}, Joyce \cite{Joyc2}, Dan Lee
\cite{Leed}, and Yng-Ing Lee \cite{Lee}. We expect the Lagrangian
self-expanders found here will also play an important role in
surgeries during Lagrangian mean curvature flow.

\subsection{Other self-similar solutions from Theorems A, B}
\label{lm33}

We now discuss the remaining solutions from Theorems A and B. As
in Remark~\ref{lm3rem1}, without loss of generality we may take
$\la_j=\pm 1$, $C=1$ and $A>0$. Section \ref{lm32} dealt with the
case $\la_1=\cdots=\la_n=1$ and $\al\ge 0$. There remain the cases
(a) $\la_1=\cdots=\la_n=1$ and $\al<0$, and (b) at least one
$\la_j$ is $-1$. In (b), we reorder $w_1,\ldots,w_n$ if necessary
so that $\la_1=\cdots=\la_m=1$ and\/ $\la_{m+1}=\cdots=\la_n=-1$.
We exclude $m=n$, as this is covered by \S\ref{lm32} and case (a),
and we exclude $m=0$, as then $L=\emptyset$. So we may take $1\le
m<n$, and the following theorem covers all the remaining cases.
\medskip

{\bf Theorem E.} {\it In Theorems A and B, suppose that either:
\begin{itemize}
\setlength{\itemsep}{0pt}
\setlength{\parsep}{0pt}
\item[{\rm(a)}]  $\la_1=\cdots=\la_n=C=1,$ $\al<0$ and\/ $A>0;$ or
\item[{\rm(b)}] $\la_1=\cdots=\la_m=1$ and\/
$\la_{m+1}=\cdots=\la_n=-1$ for some $1\le m<n,$ $C=1,$ $A>0,$
and\/~$\al\in\R$.
\end{itemize}
Then solutions exist for all\/ $s\in\R,$ and we take\/ $I=\R$. In
each of cases {\rm(a),(b)} we divide into two subcases:
\begin{itemize}
\setlength{\itemsep}{0pt}
\setlength{\parsep}{0pt}
\item[{\rm(i)}] $\sum_{j=1}^n\frac{\la_j}{\al_j}+\al=0$ and\/
$\al_1\cdots\al_n=A^2;$ or
\item[{\rm(ii)}] otherwise.
\end{itemize}
In case {\rm(i),} we have explicit solutions to \eq{lm3eq4} and
obtain \/
\begin{equation}
\begin{split}
L=\bigl\{\bigl(x_1\sqrt{\al_1}\,&e^{i(\psi_1-\la_1As/\al_1)},
\ldots,x_n\sqrt{\al_n}\,e^{i(\psi_n-\la_nAs/\al_n)}\bigr):\\
&\text{$x_1,\ldots,x_n\in\R,$ $s\in\R,$
$\ts\sum_{j=1}^n\la_jx_j^2=1$}\bigr\},
\end{split}
\label{lm3eqH}
\end{equation}
which is Hamiltonian stationary in addition to being self-similar,
and invariant under a subgroup $\R$ or $\U(1)$ of diagonal matrices
$\bigl\{\diag(e^{i\la_1t/\al_1},\allowbreak \ldots,\allowbreak
e^{i\la_nt/\al_n}):t\in\R\bigr\}$ in~$\U(n)$.

In case {\rm(ii),} $u$ and\/ $\phi-\th$ are periodic in $s$ with
period\/ $S>0$, and
\begin{align*}
u(s+S)&=u(s),& \phi_j(s+S)&=\phi_j(s)+\ga_j,\\
\phi(s+S)&=\phi(s)+\ts\sum_{j=1}^n\ga_j, &
\th(s+S)&=\th(s)+\ts\sum_{j=1}^n\ga_j,
\end{align*}
for some $\ga_1,\ldots,\ga_n\in\R$ and all\/ $s\in\R$. In case
{\rm(b)} with\/ $\al=0$ we have $\th(s)\equiv\th(0)$
and\/~$\sum_{j=1}^n\ga_j=0$.}
\medskip

The Hamiltonian stationary self-similar solutions in \eq{lm3eqH}
were obtained and studied by Lee and Wang in \cite{LeWa2}. If we
require the Lagrangian self-similar solutions in \eq{lm3eq5} to be
Hamiltonian stationary, then they must be of the form \eq{lm3eqH}
by some simple arguments.

The solutions $w_j$ obtained in Theorem E are bounded and periodic
or quasi-periodic. The {\it periodic\/} ones are much more
interesting, as then $L$ is compact in case (a), and closed in
case (b). Our next result explores this periodicity, and shows
there are many periodic solutions.
\medskip

{\bf Theorem F.} {\it In Theorem E, we say that\/ $(w_1,\ldots,w_n)$
is periodic if there exists $T>0$ with\/ $w_j(s)=w_j(s+T)$ for all\/
$s\in\R$ and\/~$j=1,\ldots,n$.

If\/ $(w_1,\ldots,w_n)$ is periodic then in case {\rm(a),} $L$ is a
compact, immersed Lagrangian self-shrinker diffeomorphic to ${\mathcal
S}^1\t{\mathcal S}^{n-1},$ and in case {\rm(b),} $L$ is a closed,
noncompact, immersed Lagrangian diffeomorphic to ${\mathcal S}^1\t{\mathcal
S}^{m-1}\t\R^{n-m},$ a self-expander if\/ $\al>0,$ a self-shrinker
if\/ $\al<0,$ and special Lagrangian if\/~$\al=0$.

In case {\rm(i),} $(w_1,\ldots,w_n)$ is periodic if and only if\/
$\frac{\la_j}{\al_j}=\mu q_j$ with\/ $\mu>0$ and\/ $q_j\in\Q$ for
$j=1,\ldots,n$. In case {\rm(ii),} $(w_1,\ldots,w_n)$ is periodic if
and only if\/ $\ga_j\in\pi\Q$ for $j=1,\ldots,n$. In both cases, for
fixed\/ $m,\al,$ there is a dense subset of initial data for which\/
$(w_1,\ldots,w_n)$ is periodic.}
\medskip

If\/ $(w_1,\ldots,w_n)$ is not periodic, then $L$ is a noncompact,
immersed Lagrangian diffeomorphic to ${\mathcal S}^{n-1}\t\R$ in case
{\rm(a)} and to ${\mathcal S}^{m-1}\t\R^{n-m+1}$ in case {\rm(b)}. It
is not closed in $\C^n,$ and the closure $\bar L$ of\/ $L$ in
$\C^n$ has dimension greater than\/~$n$.

One can  use the solutions obtained in Theorem F to form eternal
solutions of Brakke flow without mass loss, which  generalize some
of Lee and Wang's earlier results \cite{LeWa1,LeWa2}. Recall that
{\it Brakke flow} \cite{Brak} is a generalization of mean
curvature flow to {\it varifolds}, measure-theoretic
generalizations of submanifolds which may be singular. And an
eternal solution is a solution which is defined for all $t$.

For $t\in\R$, define
\begin{align*}
L_t=\,&\bigl\{(x_1\sqrt{\al_1+u(s)}\,e^{i\phi_1(s)},\ldots,
x_m\sqrt{\al_m+u(s)}\,e^{i\phi_m(s)},\\
&x_{m+1}\sqrt{\al_{m+1}-u(s)}\,e^{i\phi_{m+1}(s)},\ldots,
x_n\sqrt{\al_n-u(s)}\,e^{i\phi_n(s)}\bigr):\\
&\text{$x_1,\ldots,x_n\in\R,$ $s\in\R/T\Z$, $x_1^2+\cdots
+x_m^2-x_{m+1}^2-\cdots-x_n^2=t$}\bigr\},
\end{align*}
where $u,\phi_1,\ldots,\phi_n$ are periodic with period $T>0$. Then
$L_t$ is a closed, nonsingular, immersed Lagrangian self-expander in
$\C^n$ diffeomorphic to ${\mathcal S}^1\t{\mathcal
S}^{m-1}\t\R^{n-m}$ when $t>0,$ and a closed, nonsingular, immersed
Lagrangian self-shrinker in $\C^n$ diffeomorphic to ${\mathcal
S}^1\t{\mathcal S}^{n-m-1}\t\R^m$ when $t<0,$ and\/ $L_0$ is a
closed, immersed Lagrangian cone in $\C^n$ with link\/ ${\mathcal
S}^1\t{\mathcal S}^{m-1}\t{\mathcal S}^{n-m-1},$ with an isolated
singular point at\/~$0$.

The fact that $L_t$ form an eternal solution of Brakke flow without
mass loss is proved in Lee and Wang's paper \cite{LeWa2}.

\subsection{Translating solutions}
\label{lm34}

In \S\ref{lm31}--\S\ref{lm33} we have considered only centred
quadrics centred at 0, and only Lagrangian self-expanders and
self-shrinkers. It is an obvious question whether we can generalize
the constructions to {\it non-centred\/} quadrics on the one hand,
and to {\it Lagrangian translating solitons\/} on the other. In fact
it seems to be natural to put these ideas together, and to construct
Lagrangian translating solitons using non-centred quadrics whose
favoured axis is the direction of translation. Here is the class of
$n$-submanifolds of $\C^n$ amongst which we will seek Lagrangian
translating solitons.

\begin{ans} Fix $n\ge 2$. Consider $n$-submanifolds $L$ in $\C^n$ of
the form:
\begin{equation}
\begin{split}
L=\bigl\{\bigl(x_1w_1(s),&\ldots,x_{n-1}w_{n-1}(s),x_n+\be(s)\bigr):\\
&\text{$s\in I$, $x_1,\ldots,x_n\in\R$,
$\ts\sum_{j=1}^{n-1}\la_jx_j^2+2x_n=0$}\bigl\},
\end{split}
\label{lm3eq9}
\end{equation}
where $\la_1,\ldots,\la_{n-1}\in\R\sm\{0\}$ are nonzero constants,
$I$ is an open interval in $\R$, and $w_1,\ldots,w_{n-1}:I\ra\C
\sm\{0\}$, $\be:I\ra\C$ are smooth functions. We want $L$ to
satisfy:
\begin{itemize}
\setlength{\itemsep}{0pt}
\setlength{\parsep}{0pt}
\item[(i)] $L$ is Lagrangian;
\item[(ii)] the Lagrangian angle $\th:L\ra\R$ or
$\th:L\ra\R/2\pi\Z$ of $L$ is a function only of $s$, not of
$x_1,\ldots,x_n$; and
\item[(iii)] $L$ is a translating soliton under mean curvature flow
in $\C^n$, with translating vector $(0,\ldots,0,\al)\in\C^n$,
for~$\al\in\R$.
\end{itemize}
\label{lm3ans2}
\end{ans}

One motivation for this is the special Lagrangian submanifolds found
by the first author \cite[\S 7]{Joyc1}, which are of the form
\eq{lm3eq9}. Another is the limiting argument used to prove Theorem
G below, which recovers Ansatz \ref{lm3ans2} as a limit of Ansatz
\ref{lm3ans1}. Here is the analogue of Theorem A for this ansatz.
\medskip

{\bf Theorem G.} {\it Let\/ $\la_1,\ldots,\la_{n-1}\in\R \sm\{0\}$
and\/ $\al\in\R$ be constants, $I$ be an open interval in $\R,$
and\/ $\th:I\ra\R$ or $\th:I\ra\R/2\pi\Z,$
$w_1,\ldots,w_{n-1}:I\ra\C\sm\{0\}$ and\/ $\be:I\ra\C$ be smooth
functions. Suppose that
\begin{equation}
\begin{cases}
\begin{aligned}
\frac{\d w_j}{\d s}&=\la_je^{i\th(s)}\,\ov{w_1\cdots w_{j-1}
w_{j+1}\cdots w_{n-1}},\quad j=1,\ldots,n-1,\\
\frac{\d\th}{\d s}&=\al\Im(e^{-i\th}w_1\cdots w_{n-1}),\\
\frac{\d\be}{\d s}&=e^{i\th(s)}\,\ov{w_1\cdots w_{n-1}},
\end{aligned}
\end{cases}
\label{lm3eq10}
\end{equation}
hold in $I$. Then the submanifold\/ $L$ in $\C^n$ given by
\begin{equation}
\begin{split}
L=\bigl\{(x_1w_1(s),\ldots,x_{n-1}w_{n-1}(s),
-\ha\ts\sum_{j=1}^{n-1}\la_jx_j^2+\be(s)):&\\
x_1,\ldots,x_{n-1}\in\R,\;\> s\in I&\bigr\}
\end{split}
\label{lm3eq11}
\end{equation}
is an embedded Lagrangian diffeomorphic to $\R^n,$ with Lagrangian
angle $\th(s)$. When $\al\ne 0,$ it is a Lagrangian translating
soliton with translating vector $(0,\ldots,0,\al)\in\C^n$. When
$\al=0,$ it is special Lagrangian, and the construction reduces to
that of Joyce~{\rm\cite[\S 7]{Joyc1}}.}
\medskip

We can prove this directly, following the proof of Theorem A in
\S\ref{lm4}. This is straightforward, and we leave it as an exercise
for the interested reader. Instead, we give a somewhat informal
proof which obtains Theorem G from Theorem A by a limiting
procedure, since this gives more insight into why the construction
should generalize in this way.

\begin{proof}[Proof of Theorem G, assuming Theorem A] Let
$\la_1,\ldots,\la_{n-1}\in\R\sm\{0\}$ and $\al\in\R$ be constants,
$I$ be an open interval in $\R$, and $\th:I\ra\R$ or
$\th:I\ra\R/2\pi\Z,$ $w_1,\ldots,w_{n-1}:I\ra\C\sm\{0\}$ and
$\be:I\ra\C$ be smooth functions. Let $R>0$. Define constants
$\ti\la_1,\ldots,\ti\la_n,\ti C\in\R\sm\{0\}$ and $\ti\al\in\R$, an
open interval $\ti I$, and smooth $\ti\th:\ti I\ra\R$ or $\ti\th:\ti
I\ra\R/2\pi\Z,$ $\ti w_1,\ldots,\ti w_n:\ti I\ra\C\sm\{0\}$ by
\begin{equation}
\begin{gathered}
\ti\la_j=\la_j,\;\> j=1,\ldots,n-1,\;\> \ti\la_n=R,\;\> \ti
C=R,\;\> \ti\al=\al,\;\> \ti I=R^{-1}I,\\
\ti s=R^{-1}s,\;\>
\ti w_j(\ti s)=w_j(R\ti s)=w_j(s),\;\> j=1,\ldots,n-1,\\
\ti w_n(\ti s)=R+\be(R\ti s)=R+\be(s),\;\> \ti\th(\ti s)=\th(R\ti
s)=\th(s).
\end{gathered}
\label{lm3eq12}
\end{equation}
We suppose that $\be\ne -R$ so that $\ti w_n$ maps~$\ti
I\ra\C\sm\{0\}$.

Apply Theorem A to this new data $\ti\la_j,\ti C,\ti\al,\ti
I,\ti\th,\ti w_j$. This yields o.d.e.s \eq{lm3eq2} upon $\ti
w_j,\ti\th$, in terms of derivatives with respect to $\ti s$, and
defines \eq{lm3eq3} a self-similar Lagrangian $\ti L$ in $\C^n$ when
these o.d.e.s hold. Define $L=\ti L-(0,\ldots,0,R)$, that is, $L$ is
$\ti L$ translated by the vector $-(0,\ldots,0,R)$. Rewriting the
o.d.e.s \eq{lm3eq2} in terms of $R,\la_j,\al,I,\th,w_j,\be$ using
\eq{lm3eq12} yields
\begin{equation}
\begin{cases}
\frac{\d w_j}{\d s}&\!\!\!\!=\la_je^{i\th(s)}\,\ov{w_1\cdots w_{j-1}
w_{j+1}\cdots w_{n-1}}\,(1+R^{-1}\ov\be(s)),
\>\;\text{all $j$,}\\
\frac{\d\th}{\d s}&\!\!\!\!=\al\Im\bigl(e^{-i\th}w_1\cdots w_{n-1}
(1+R^{-1}\be(s))\bigr),\\
\frac{\d\be}{\d s}&\!\!\!\!=e^{i\th(s)}\,\ov{w_1\cdots w_{n-1}}.
\end{cases}
\label{lm3eq13}
\end{equation}
Rewriting $\ti L$ in \eq{lm3eq3} in terms of
$R,\la_j,\al,I,\th,w_j,\be$, translating by $-(0,\ab\ldots,\ab 0,R)$
to get $L$, and replacing $x_n$ in \eq{lm3eq3} by
$1+R^{-1}\bar{x}_n$, yields
\begin{equation}
\begin{split}
L=&\bigl\{\bigl((x_1w_1(s),\ldots,x_{n-1}w_{n-1}(s),
\bar{x}_n+\be(s)+R^{-1}\be(s)\bar{x}_n\bigr):\\
&\text{$x_1,\ldots,x_{n-1},\bar{x}_n\!\in\!\R$, $s\!\in\! I$,
$\ts\sum_{j=1}^{n-1}\la_jx_j^2\!+\!2\bar{x}_n+R^{-1}\bar{x}_n^2
\!=\!0$}\bigr\}.
\end{split}
\label{lm3eq14}
\end{equation}

The conclusion of Theorem A is that $\ti L$ satisfies $\al
F^\perp=RH$. Since $L$ is the translation of $\ti L$ by
$-(0,\ldots,0,R)$, and this subtracts $(0,\ldots,0,R)$ from $F$, we
see that $L$ satisfies $\al\bigl(F+(0,\ldots,0,R)\bigr){}^\perp=RH$.
Dividing by $R$ and setting $T=(0,\ldots,0,\al)$, this shows that
$L$ satisfies $H=T^\perp+R^{-1}\al F^\perp$.

Now let us take the limit $R\ra\iy$. Then \eq{lm3eq13} reduces to
\eq{lm3eq10}, as the $R^{-1}$ terms disappear, and \eq{lm3eq14}
reduces to \eq{lm3eq11}, as
$\bar{x}_n=-\ha\ts\sum_{j=1}^{n-1}\la_jx_j^2$ in the limit. The
equation $H=T^\perp+R^{-1}\al F^\perp$ for $L$ becomes
$H=T^\perp$, so $L$ is a Lagrangian translating soliton with
translating vector~$T=(0,\ldots,0,\al)$.

It remains to show that $L$ is {\it embedded}, that is, the
immersion $\io:(x_1,\ab\ldots,\ab x_{n-1},\ab s)\mapsto(x_1w_1(s),
\ldots, x_{n-1}w_{n-1}(s),-\ha\sum_{j=1}^{n-1}\la_jx_j^2+\be(s))$ is
injective. Combining \eq{lm3eq13} with equation \eq{lm3eq6} of
Theorem B gives $\Im\frac{\d\be}{\d s}\equiv -Ae^{-\al u/2}$ for
some $A\in\R$. Thus $\Im\be$ is strictly decreasing in $s$ if $A>0$,
and strictly increasing if $A<0$. In both cases, if
$\io(x_1,\ldots,x_{n-1},s)=(z_1,\ldots,z_n)$ then $\Im z_n$
determines $s$, and given $s$, we have $x_j=w_j(s)^{-1}z_j$ for
$j=1,\ldots,n-1$. So $\io$ is injective if $A\ne 0$. When $A=0$ we
can solve explicitly and show $\io$ is injective.
\end{proof}

Note that the first two equations in \eq{lm3eq10} are exactly the
same as \eq{lm3eq2}, replacing $n$ by $n-1$. Having chosen some
solutions $w_1,\ldots,w_{n-1},\th$ to the first two equations of
\eq{lm3eq10}, the third equation of \eq{lm3eq10} determines $\be$
uniquely up to $\be\mapsto\be+c$, by integration. Actually we can
write $\be$ explicitly in terms of $u,\th$: in the notation of
Theorem B, if $\al\ne 0$ then \eq{lm3eq4} and the last equation of
\eq{lm3eq10} give $\frac{\d\be}{\d s}=\ha\frac{\d u}{\d s}-
\frac{i}{\al}\frac{\d\th}{\d s}$, which integrates to $\be(s)=\ha
u(s)-\frac{i}{\al}\th(s)+K$, for $K\in\C$. So we deduce:
\medskip

{\bf Corollary H.} {\it In the situation of Theorem G, when $\al\ne
0,$ the Lagrangian translating soliton $L$ may be rewritten
\begin{equation}
\begin{split}
\!\!L\!=\!\bigl\{\bigl(x_1\sqrt{\al_1\!+\!\la_1u(s)}\,e^{i\phi_1(s)}\!,
\ldots,x_{n-1}\sqrt{\al_{n-1}\!+\!\la_{n-1}u(s)}\,e^{i\phi_{n-1}(s)}&,\!\! \\
\ts\ha u(s)-\ha\sum_{j=1}^{n-1}\la_jx_j^2-\frac{i}{\al}\th(s)+K
\bigr):x_1,\ldots,x_{n-1}\in\R,\;\> s\in I\bigr\}&,\!\!
\end{split}
\label{lm3eq15}
\end{equation}
where $K\in\C$ and\/ $u,\al_1,\ldots,\al_{n-1},\phi_1,\ldots,
\phi_{n-1}$ are as in Theorem B with\/ $n-1$ in place of\/ $n,$ and
satisfy \eq{lm3eq4} and\/ \eq{lm3eq6} for some\/~$A\in\R$.}
\medskip

Proposition \ref{lm2prop2} implies that the Lagrangian angle $\th$
of $L$ should be of the form $-\al\Im z_n\vert_L+c$, where
$(z_1,\ldots,z_n)$ are the complex coordinates on $\C^n$. The
imaginary part of the last coordinate in \eq{lm3eq15} is
$-\frac{1}{\al} \th(s)+\Im K$, so the proposition holds
with~$c=\al\Im K$.

Theorems C, D, E and F give a good description of solutions of
\eq{lm3eq4}, and hence of the Lagrangian translating solitons $L$
from Theorem G and Corollary H. In the authors' opinion, the most
interesting case of Theorem G is when $\la_1,\ldots,\la_{n-1}>0$
and $\al\ge 0$. (This is equivalent to the case
$\la_1,\ldots,\la_{n-1} <0$ and $\al\le 0$, changing the sign of
the last coordinate in $\C^n$.) The following result combines
Theorems C, D with $n-1$ in place of $n$, Theorem G and Corollary
H. For simplicity we set $\psi_1=\cdots=\psi_{n-1}=0$ and~$K=-\ha
u_*$, where $u_*$ is defined in the proof of Theorem~C (see
Figure~\ref{lm5fig}).
\medskip

{\bf Corollary I.} {\it For given constants $\al\ge 0$ and\/
$a_1,\ldots,a_{n-1}>0,$ define
\begin{equation*}
\phi_j(y)=\int_0^y\frac{\d t}{(\frac{1}{a_j}+t^2)\sqrt{P(t)}}\,,
\;\>\text{where}\;\>
P(t)=\frac{1}{t^2}\Big(\prod_{k=1}^{n-1}(1+a_kt^2)e^{\al
t^2}-1\Big),
\end{equation*}
for $j=1,\ldots,n-1$ and\/ $y\in\R$. Then when $\al\ne 0,$
\begin{align}
L\!=\!\bigl\{\bigl(&x_1\ts\sqrt{\frac{1}{a_1}\!+\!y^2}\,e^{i\phi_1(y)},
\ldots,x_{n-1}\sqrt{\frac{1}{a_{n-1}}\!+\!y^2}\,e^{i\phi_{n-1}(y)},
\ts\ha y^2\!-\!\ha\sum_{j=1}^{n-1}x_j^2
\nonumber\\
&-\ts\frac{i}{\al}\sum_{j=1}^{n-1}\phi_j(y)-\ts\frac{i}{\al}
\arg(y+iP(y)^{-1/2})\bigr):x_1,\ldots,x_{n-1},y\in\R\bigr\}
\label{lm3eq16}
\end{align}
is a closed, embedded Lagrangian in $\C^n$ diffeomorphic to
$\R^n,$ which is a Lagrangian translating soliton with translating
vector $(0,\ldots,0,\al)\in\C^n$. When $\al=0$,
\begin{align*}
L\!=\!\bigl\{\bigl(&x_1\ts\sqrt{\frac{1}{a_1}\!+\!y^2}\,e^{i\phi_1(y)},
\ldots,x_{n-1}\sqrt{\frac{1}{a_{n-1}}\!+\!y^2}\,e^{i\phi_{n-1}(y)},
\ts\ha y^2\!-\!\ha\sum\limits_{j=1}^{n-1}x_j^2 \nonumber\\&
+i\int_0^y \frac{\d
t}{\sqrt{\ts\frac{1}{t^2}\bigl(\prod_{k=1}^n(1+a_kt^2)-1\bigl)}}\bigr):
x_1,\ldots,x_{n-1},y\in\R\bigr\}
\end{align*}
is a closed, embedded Lagrangian in $\C^n$ diffeomorphic to $\R^n,$
which is special Lagrangian.

There exist\/ $\bar\phi_1,\ldots,\bar\phi_{n-1}\in
(0,\frac{\pi}{2}]$ such that\/ $\phi_j(y)\ra\bar\phi_j$ as $y\ra\iy$
and\/ $\phi_j(y)\ra -\bar\phi_j$ as $y\ra-\iy$ for $j=1,\ldots,n-1$.
These satisfy $\bar\phi_1+\cdots+\bar\phi_{n-1}<\frac{\pi}{2}$ if\/
$\al>0,$ and\/ $\bar\phi_1+\cdots+\bar\phi_{n-1}=\frac{\pi}{2}$ if\/
$\al=0$. For fixed $\al>0,$ the map $(a_1,\ldots,a_{n-1})\mapsto
(\bar\phi_1,\ldots,\bar\phi_{n-1})$ is a {\rm 1-1} correspondence
from $(0,\iy)^{n-1}$ to
$\bigl\{(\bar\phi_1,\ldots,\bar\phi_{n-1})\in
(0,\frac{\pi}{2})^{n-1}:\bar\phi_1+\cdots+\bar\phi_{n-1}<
\frac{\pi}{2}\bigr\}$. When $\al=0,$ the map has
image~$\bigl\{(\bar\phi_1,\ldots,\bar\phi_{n-1})\in
(0,\frac{\pi}{2}]^{n-1}:\bar\phi_1+\cdots+\bar\phi_{n-1}=
\frac{\pi}{2}\bigr\}$.

The Lagrangian angle of\/ $L$ in \eq{lm3eq16} varies between
$\sum_{j=1}^{n-1}\bar{\phi}_j$ and\/
$\pi-\sum_{j=1}^{n-1}\bar{\phi}_j$. Thus, when $\al>0,$ by choosing
$\sum_{j=1}^{n-1}\bar{\phi}_j$ close to $\frac{\pi}{2},$ the
oscillation of the Lagrangian angle of\/ $L$ can be made arbitrarily
small.}
\medskip

We can give the following heuristic description of $L$ in
\eq{lm3eq16}. If $y\gg 0$ then $\phi_j(y)\approx\bar\phi_j$ and
$\sqrt{\frac{1}{a_j}+y^2}\approx y$, and the terms $-\frac{i}{\al}
\sum_{j=1}^n\phi_j(y)-\frac{i}{\al}\arg(y+iP(y)^{-1/2})$ are
negligible compared to $\ha y^2$ in the last coordinate. Thus, the
region of $L$ with $y\!\gg\! 0$ is in a weak sense approximate~to
\begin{equation*}
\bigl\{\bigl(x_1ye^{i\bar\phi_1},\ldots,
x_{n-1}ye^{i\bar\phi_{n-1}},\ts\ha y^2-\ha\sum
\limits_{j=1}^{n-1}x_j^2\bigr):x_1,\ldots,x_{n-1}\in\R,\; y>0\bigr\}.
\end{equation*}
But this is just an unusual way of parametrizing
\begin{equation*}
L_1=\bigl\{\bigl(y_1e^{i\bar\phi_1},\ldots,
y_{n-1}e^{i\bar\phi_{n-1}},y_n):y_j\in\R\bigr\}\sm
\bigl\{(0,\ldots,0,y_n):y_n\le 0\bigr\},
\end{equation*}
the complement of a ray in a Lagrangian plane. Similarly, the region
of $L$ with $y\ll 0$ is in a weak sense approximate to
\begin{equation*}
L_2\!=\!\bigl\{\bigl(y_1e^{-i\bar\phi_1},\ldots,
y_{n-1}e^{-i\bar\phi_{n-1}},y_n):y_j\!\in\!\R\bigr\}\sm
\bigl\{(0,\ldots,0,y_n):y_n\!\le\!0\bigr\}.
\end{equation*}

So, $L$ can be roughly described as asymptotic to the union of two
Lagrangian planes $L_1,L_2\cong\R^n$ which intersect in an $\R$ in
$\C^n$, the $y_n$-axis $\bigl\{(0,\ldots,0,y_n):y_n\in\R\bigr\}$. To
make $L$, we glue these Lagrangian planes by a kind of `connect sum'
along the negative $y_n$-axis $\bigl\{(0,\ldots,0,y_n):y_n\le
0\bigr\}$. Under Lagrangian mean curvature flow, $L_1,L_2$ remain
fixed, but the gluing region translates in the positive $y_n$
direction, as though $L_1,L_2$ are being `zipped together'.

Note too that when the oscillation of the Lagrangian angle of $L$ is
small compared to $\md{\al}$, we see from Proposition \ref{lm2prop2}
that $\Im z_n\vert_L$ is confined to a small interval, where
$(z_1,\ldots,z_n)$ are the complex coordinates on $\C^n$. That is,
$L$ is close to the affine $\R^{2n-1}$ in $\C^n$ defined by~$\Im
z_n=-\frac{\pi}{2\al}$.

\begin{quest} Can the translating solitons with small Lagrangian
angle oscillation in Corollary I arise as blow-ups of finite time
singularities for Lagrangian mean curvature flow, particularly
when~$n=2$?
\end{quest}

It is important to answer this question in developing a regularity
theory for the flow. Such relations have been observed before in
codimension one mean curvature by White \cite{Whit,Whit1}, and
Huisken and Sinestrari \cite{HuSi}, and in Ricci flow by Perelman
\cite{Pere}. See also the recent work by Neves and Tian \cite{Neve1}
for related discussions.

We can also ask about the Lagrangian translating solitons from
Theorem G and Corollary H coming from Theorem E with $n-1$ in
place of $n$. As in Theorem~E, we take $I=\R$. Using the notation
of Theorem B for $w_1,\ldots,w_{n-1}, \th$, observe that the third
equation of \eq{lm3eq10} gives
\begin{equation}
\Im\frac{\d\be}{\d s}=-Q(u)^{1/2}\sin(\phi-\th)=-Ae^{-\al u/2}.
\label{lm3eq17}
\end{equation}

As in Remark \ref{lm3rem1}(d), when $A=0$ the Lagrangian $L$ is an
open subset of an affine Lagrangian plane $\R^n$ in $\C^n$, which is
not interesting, so we restrict to $A\ne 0$. Then \eq{lm3eq17} shows
that either $\Im\frac{\d\be}{\d s}>0$ for all $s\in\R$, or
$\Im\frac{\d\be}{\d s}<0$ for all $s\in\R$. Thus $\be$ can never be
periodic, so we have no analogue of Theorem F in the translating
case. We can also deduce from this that the Lagrangians are closed,
embedded, diffeomorphic to $\R^n$, and when $\al\ne 0$ have {\it
infinite oscillation of the Lagrangian angle}. This implies that
these Lagrangian translating solitons from Theorems G and E {\it
cannot arise as blow-ups of finite time singularities for Lagrangian
mean curvature flow}.

\section{A construction for self-similar Lagrangians}
\label{lm4}

We now prove:
\smallskip

{\bf Theorem A.} {\it Let\/ $\la_1,\ldots,\la_n,C\in\R\sm\{0\}$
and\/ $\al\in\R$ be constants, $I$ be an open interval in $\R,$
and\/ $\th:I\ra\R$ or $\th:I\ra\R/2\pi\Z$ and\/
$w_1,\ldots,w_n:I\ra\C\sm\{0\}$ be smooth functions. Suppose that
\begin{equation}
\begin{aligned}
\frac{\d w_j}{\d s}&= \la_je^{i\th(s)}\,\ov{w_1\cdots w_{j-1}
w_{j+1}\cdots w_n},\qquad j=1,\ldots,n,\\
\frac{\d\th}{\d s}&=\al\Im(e^{-i \th(s)} w_1\cdots w_n),
\end{aligned}
\label{lm4eq1}
\end{equation}
hold in $I$. Then the submanifold\/ $L$ in $\C^n$ given by
\begin{equation}
L=\bigl\{\bigl(x_1w_1(s),\ldots,x_nw_n(s)\bigr):\text{$s\in I,$
$x_j\in\R,$ $\ts\sum_{j=1}^{n}\la_jx_j^2=C$}\bigl\},
\label{lm4eq2}
\end{equation}
is Lagrangian, with Lagrangian angle $\th(s)$ at\/
$(x_1w_1(s),\ldots,x_nw_n(s)),$ and its position vector\/ $F$ and
mean curvature vector \/ $H$ satisfy\/ $\al F^\perp=CH$. That is,
$L$ is a self-expander when\/ $\al/C>0$ and a self-shrinker when\/
$\al/C<0$. When $\al=0$ the Lagrangian angle $\th$ is constant, so
that\/ $L$ is special Lagrangian, with\/ $H=0$. In this case the
construction reduces to that of Joyce~{\rm\cite[\S 5]{Joyc1}}.}
\smallskip

\begin{proof} Define $\Si=\bigl\{(x_1,\ldots,x_n)\in\R^n:
\sum_{j=1}^n\la_jx_j^2=C\bigr\}$. Then $\Si$ is a nonsingular quadric in
$\R^n$, an $(n-1)$-manifold. Define a smooth map $\io:\Si\t
I\ra\C^n$ by $\io:\bigl((x_1,\ldots,x_n),s\bigr)\longmapsto
\bigl(x_1w_1(s),\ldots,x_nw_n(s)\bigr)$. Then $L=\io(\Si\t I)$. The
determinant calculation below implies $\io$ is an {\it immersion},
and so $L$ is a {\it nonsingular immersed $n$-submanifold\/}
in~$\C^n$.

Fix $\bs x=(x_1,\ldots,x_n)\in\Si$ and $s\in I$. We will find the
tangent space $T_{\io(\bs x,s)}L$, show that it is Lagrangian, and
compute its Lagrangian angle. Let $e_1,\ldots,e_{n-1}$ be an
orthonormal basis for $T_{\bs x}\Si$ in $\R^n$, and write
$e_j=(a_{j1},\ldots,a_{jn})$ in $\R^n$ for $j=1,\ldots,n-1$. Let
$e_n=\bigl(\sum_{j=1}^{n}\la_j^2 x_j^2\bigr){}^{-1/2}\cdot\ab
(\la_1x_1,\ab \ldots,\ab\la_nx_n)$. Then $e_n$ is a unit normal
vector to $\Si$ at $\bs x$ in $\R^n$. Let $e_1,\ldots,e_{n-1}$ be
chosen so that $e_1,\ldots,e_{n-1},e_n$ is an oriented orthonormal
basis for $\R^n$. Then $\det(e_1\,\ldots\,e_n)=1$, regarding
$e_1,\ldots,e_n$ as column vectors, and $(e_1\,\ldots\,e_n)$ as an
$n\t n$ matrix.

Now $e_1,\ldots,e_{n-1},\frac{\pd}{\pd s}$ is a basis for $T_{(\bs
x,s)}(\Si\t I)$. Define $f_j=\d\io(e_j)\in\C^n$ for $j=1,\ldots,n-1$
and $f_n=\d\io(\frac{\pd}{\pd s})\in\C^n$. Then $f_1,\ldots,f_n$ is
a basis for $T_{\io(\bs x,s)}L$, over $\R$. From the definitions we
have $f_j=\bigl(a_{j1}w _1(s),\ldots,a_{jn}w _n(s)\bigr)$ for $j<n$,
and $f_n=\bigl(x_1\dot{w}_1(s),\ldots, x_n\dot{w}_n (s)\bigr)$.
Therefore
\begin{align*}
\an{f_j,Jf_k}&= \Re\bigl(-i\ts\sum_{l=1}^{n}a_{jl}a_{kl}\ms{w
_l}\bigr)=0, \quad \text{for $j,k=1,\ldots,n-1$}, \\
\an{f_j,Jf_n}&=\Re\bigl(-i\,w_{1}\cdots w_{n}e^{-i\th}
\ts\sum_{j=1}^{n}a_{jl}\la_lx_l)=0, \;\> \text{for
$j\!=\!1,\ldots,n\!-\!1$},
\end{align*}
where in the second equation we use the first equation of
\eq{lm4eq1} to substitute for $\dot{w}_l(s)$, and the fact that
$(\la_1x_1,\ldots,\la_nx_n)$ is normal to $\Si$ at $\bs x$, and so
orthogonal to~$e_j=(a_{j1},\ldots,a_{jn})$.

Thus $\an{f_j,Jf_k}=0$ for $j,k=1,\ldots,n$, so the symplectic form
$\omega(*,*)=\an{*,J*}$ on $\C^n$ vanishes on
$\an{f_1,\ldots,f_n}_\R=T_{\io(\bs x,s)}L$, and $T_{\io(\bs x,s)}L$
is a Lagrangian plane in $\C^n$. Hence $L$ is {\it Lagrangian}. To
compute the Lagrangian angle, write $w _j(s)=r_{j}(s)e^{i\phi
_j(s)}$ and $\phi(s)=\sum_{j=1}^{n} \phi _j(s)$, where $r_j(s)=\md{w
_j(s)}$. Then
\begin{small}
\begin{align*}
&\det(f_1\,\cdots\,f_n) \\
&=\left\vert\begin{matrix} a_{11}r_1e^{i\phi_1(s)} & \cdots &
a_{(n-1)1}r_1e^{i\phi_1(s)} &
x_1\la_1 r_2\cdots r_ne^{i(\phi_1(s)+\th(s)-\phi(s))} \\
a_{12}r_2e^{i\phi_2(s)} & \cdots & a_{(n-1)2}r_2e^{i\phi_2(s)} &
x_2\la_2 r_1r_3\cdots r_ne^{i(\phi_2(s)+\th(s)-\phi(s))} \\
\vdots & \vdots & \vdots & \vdots \\
a_{1n}r_ne^{i\phi_n(s)} & \cdots & a_{(n-1)n}r_ne^{i\phi_n(s)} &
x_n\la_n r_1\cdots r_{n-1}e^{i(\phi_n(s)+\th(s)-\phi(s))}
\end{matrix}\right\vert\\
&= r_1^2\cdots r_n^2e^{i\th}\left\vert\begin{matrix}
a_{11} & \cdots & a_{(n-1)1} & \frac{x_1\la_1}{r_1^2}\\
a_{12} & \cdots & a_{(n-1)2} & \frac{x_2\la_2}{r_2^2}\\
\vdots & \vdots & \vdots & \vdots \\
a_{1n} & \cdots & a_{(n-1)n} & \frac{x_n\la_n}{r_n^2}
\end{matrix}\right\vert
=\frac{r_1^2\cdots r_n^2e^{i\th}}{\sqrt{\sum_{l=1}^{n}\la_l^2
x_l^2}}\sum_{l=1}^{n}\frac{\la_l^2x_l^2}{r_l^2}.
\end{align*}
\end{small}

Here in the second step we extract factors of $r_je^{i\phi_j(s)}$
from the $j^{\rm th}$ row for $j=1,\ldots,n$, and a factor $
r_1\cdots r_ne^{i(\th(s)-\phi(s))}$ from the $n^{\rm th}$ column.
The factors $e^{i\phi_1(s)}e^{i\phi_2(s)}\cdots
e^{i\phi_n(s)}e^{-i\phi(s)}$ cancel to give 1. In the third and
final step, we note that the first $n-1$ columns of the matrix on
the third line are $e_1,\ldots,e_{n-1}$, and $e_1,\ldots,e_n$ are
orthonormal with $\det(e_1\,\ldots\,e_n)=1$, so we project the
vector $(\frac{\la_1x_1}{r_1^2},\ldots,\frac{\la_nx_n}{r_n^2})$ to
$e_n$ to calculate the determinant. This shows that the Lagrangian
angle on $L$ at $(x_1w_1(s),\ldots,x_nw_n(s))$ is $\th(s)$, as we
have to prove. Also, as it shows that $\det(f_1\,\cdots\,f_n)\ne 0$,
this calculation implies that $\d\io:T_{(\bs x,s)}(\Si\t I)\ra\C^n$
is injective, and $\io$ is an {\it immersion}, as we claimed above.

The matrix $(g_{ab})$ of the metric on $L$ w.r.t.\ the basis
$f_1,\ldots,f_n$ is
\begin{equation}
\ts g_{nn}\!=\!r_1^2\cdots r_n^2\, \sum_{l=1}^{n}\frac{\la_l^2
x_l^2}{r_l^2}, \; g_{jn}\!=\!g_{nj}\!=\!0, \; \text{and}\;
g_{jk}\!=\!\sum_{l=1}^{n}a_{jl}a_{kl}r_l^{2}
\label{lm4eq3}
\end{equation}
for $j,k=1,\ldots,n-1$. Hence by \eq{lm4eq1} the {\it mean curvature
vector\/} is
\begin{equation}
H=J\nabla\th=\frac{\dot{\th}}{g_{nn}}\,Jf_n=\frac{\al r_1\cdots
r_{n}\sin(\phi-\th)}{g_{nn}}\,Jf_n.
\label{lm4eq4}
\end{equation}

The normal projection of the position vector $F$ is computed by
\begin{align*}
\an{F,Jf_l}&=\Re\bigl(-i\ts\sum_{j=1}^{n}r_{j}^{2}x_ja_{lj}\bigr)=0,\\
\an{F,Jf_n}&=\Re\bigl(-ir_{1}\cdots
r_{n}e^{i(\phi-\th)}\ts\sum_{j=1}^{n} \la_j x_j^2\bigr)=Cr_{1}\cdots
r_{n}\sin(\phi-\th).
\end{align*}
It follows that
\begin{equation}
F^\perp=\frac{Cr_{1}\cdots r_{n}\sin(\phi-\th)}{g_{nn}}\,Jf_n.
\label{lm4eq5}
\end{equation}
Equations \eq{lm4eq4} and \eq{lm4eq5} give $\al F^\perp=CH$, as we
have to prove.
\end{proof}

We can rewrite, simplify, and partially solve the
equations~\eq{lm4eq1}.
\medskip

{\bf Theorem B.} {\it In the situation of Theorem A, let\/
$w_1,\ldots,w_n,\th$ satisfy\/ \eq{lm4eq1}. Write $w_j\equiv
r_je^{i\phi_j}$ and\/ $\phi=\sum_{j=1}^n\phi_j,$ for functions
$r_j:I\ra(0,\iy)$ and\/ $\phi_1,\ldots,\phi_n,\phi:I\ra\R$ or
$\R/2\pi\Z$. Fix $s_0\in I$. Define $u:I\ra\R$ by
\begin{equation}
u(s)=2\int_{s_0}^sr_1(t)\cdots
r_n(t)\cos\bigl(\phi(t)-\th(t)\bigr)\d t.
\label{lm4eq6}
\end{equation}
Then $r_j^2(s)\equiv\al_j+\la_ju(s)$ for $j=1,\ldots,n$ and\/ $s\in
I,$ where $\al_j=r_j^2(s_0)$. Define a degree $n$ polynomial $Q(u)$
by $Q(u)=\prod_{j=1}^n(\al_j+\la_ju)$. Then the system of equations
\eq{lm4eq1} can be rewritten as
\begin{equation}
\begin{cases}
\begin{aligned}
\frac{\d u}{\d s}&=2Q(u)^{1/2}\cos(\phi-\th),\\
\frac{\d\phi_j}{\d s} &= -\frac{\la_jQ(u)^{1/2}\sin(\phi-\th)}{
\al_j+ \la_ju}\,, \qquad j=1,\ldots,n,\\
\frac{\d\phi}{\d s}&=-Q(u)^{1/2}(\ln Q(u))' \sin(\phi-\th),\\
\frac{\d\th}{\d s} &= \al Q(u)^{1/2}\sin(\phi-\th).
\end{aligned}
\end{cases}
\label{lm4eq7}
\end{equation}
The Lagrangian self-similar solution $L$ in Theorem A may be
rewritten
\begin{equation}
\begin{split}
L=\bigl\{(x_1\sqrt{\al_1+\la_1u(s)}\,&e^{i\phi_1(s)},
\ldots,x_n\sqrt{\al_n+\la_nu(s)}\,e^{i\phi_n(s)}):\\
&\text{$x_1,\ldots,x_n\in\R,$ $s\in I,$
$\ts\sum_{j=1}^n\la_jx_j^2=C$}\bigr\}.
\end{split}
\label{lm4eq8}
\end{equation}
Moreover, for some $A\in\R$ the equations \eq{lm4eq7} have the first
integral}
\begin{equation}
Q(u)^{1/2}e^{\al u/2}\sin(\phi-\th)\equiv A.
\label{lm4eq9}
\end{equation}
\smallskip

\begin{proof} Using equations \eq{lm4eq1} and \eq{lm4eq7}, for
$j=1,\ldots,n$ we have
\begin{align*}
\frac{\d(r_j^2)}{\d s}&=\frac{\d(\ms{w_j})}{\d
s}\!=\!w_j\frac{\d\bar w_j}{\d s}\!+\!\bar w_j\frac{\d
w_j}{\d s}\!=\!\la_je^{-i\th}w_1\cdots w_n\!+\!
\la_je^{i\th}\ov{w_1\cdots w_n}\\
&=2\la_j \Re(e^{i(\phi-\th)}r_1\cdots r_n)=
2\la_j\cos(\phi-\th)\,r_1\cdots r_n=\la_j\frac{\d u}{\d s}\,.
\end{align*}
Thus $r_j^2-\la_ju$ is constant in $I$, and at $s=s_0$ we have
$r_j^2(s_0)=\al_j$ and $u(s_0)=0$, so $r_j^2(s)\equiv\al_j+
\la_ju(s)$ for $j=1,\ldots,n$ and $s\in I$, as we have to prove.

Differentiating \eq{lm4eq6} gives $\frac{\d u}{\d s}=2r_1\cdots
r_n\cos(\phi-\th)$. But $r_j^2=\al_j+\la_ju$ and the definition of
$Q$ imply that $Q(u)=\prod_{j=1}^nr_j^2$, so $\frac{\d u}{\d s}=2
Q(u)^{1/2}\cos(\phi-\th)$, the first equation of \eq{lm4eq7}. As
$w_j=r_je^{i\phi_j}$ we have
\begin{equation*}
\frac{\d w_j}{\d s}=\frac{\d r_j}{\d
s}\,e^{i\phi_j}+ir_je^{i\phi_j}\frac{\d\phi_j}{\d s}.
\end{equation*}
Thus $r_j\frac{\d\phi_j}{\d s}=\Im(e^{-i\phi _j}\frac{\d w _j}{\d
s})$, and the second equation of \eq{lm4eq7} follows from the first
equation of \eq{lm4eq1}, $w_j=r_je^{i\phi_j}$,
$\phi=\sum_{j=1}^n\phi_j$ and $Q(u)^{1/2}=r_1\cdots r_n$.

Summing the second equation of \eq{lm4eq7} over $j=1,\ldots,n$ gives
\begin{equation*}
\frac{\d\phi}{\d s}=- Q(u)^{1/2}\sum_{j=1}^n
\frac{\la_j\sin(\phi-\th)}{\al_j+\la_ju}=-Q(u)^{1/2}(\ln
Q(u))'\sin(\phi-\th),
\end{equation*}
the third equation of \eq{lm4eq7}. The final equation of
\eq{lm4eq7} follows from the second equation of \eq{lm4eq2},
$w_j=r_je^{i\phi_j}$, $\phi=\sum_{j=1}^n\phi_j$ and
$Q(u)^{1/2}=r_1\cdots r_n$. Equation~\eq{lm4eq8} is immediate from
\eq{lm4eq2} and $w_j=r_je^{i\phi_j}$. Finally, using \eq{lm4eq7}
we find that
\begin{align*}
&\frac{\d}{\d s}\bigl(Q(u)^{1/2}e^{\al u/2}\sin(\phi-\th)\bigr)= \ha
Q(u)^{1/2}\ln(Q(u))'\frac{\d u}{\d
s}\,e^{\al u/2}\sin(\phi-\th)\\
&+Q(u)^{1/2}\frac{\al}{2}e^{\al u/2}\frac{\d u}{\d s} \sin(\phi-\th)
+Q(u)^{1/2}e^{\al u/2}\cos(\phi-\th)\frac{\d(\phi-\th)}{\d s}\\
&=Q(u)^{1/2}e^{\al u/2}\bigl[\sin(\phi-\th)
\bigl(\ha\ln(Q(u))'+\ts\frac{\al}{2}\bigr)2Q(u)^{1/2}\cos(\phi-\th)\\
&+\!\cos(\phi\!-\!\th)\bigl(-Q(u)^{1/2}(\ln Q(u))' \sin(\phi\!-\!\th)
\!-\!\al Q(u)^{1/2}\sin(\phi\!-\!\th)\bigr)\bigr]\!=\!0.
\end{align*}
So the left hand side of \eq{lm4eq9} is a constant, say $A$ in $\R$.
\end{proof}

\section{Self-expanders diffeomorphic to ${\mathcal S}^{n-1}\t \R$}
\label{lm5}

We now prove Theorems C and D of~\S\ref{lm32}.
\medskip

{\bf Theorem C.} {\it In Theorems A and B, suppose that\/
$\la_1=\cdots=\la_n=C=1$, $\al\ge 0$ and\/ $A<0$. Then any
solution of\/ {\rm\eq{lm4eq1},} or equivalently of\/
{\rm\eq{lm4eq7},} on an interval\/ $I$ in $\R$ can be extended to
a unique largest open interval\/ $I_{\max}$ in $\R$. Take
$I=I_{\max}$. Then by changing variables from $s$ in $I_{\max}$ to
$y=y(s)$ in $\R,$ we may rewrite the Lagrangian self-expander $L$
of\/ \eq{lm4eq2} and\/ \eq{lm4eq8} explicitly as follows.
Conversely, every $L$ of the following form comes from Theorems~A
and~B with\/ $\la_1=\cdots=\la_n=C=1$, $\al\ge 0$ and\/~$A<0$.

For given constants $\al\ge 0,$ $a_1,\ldots,a_n>0$ and\/
$\psi_1,\ldots,\psi_n\in\R,$ define $w_j(y)=e^{i\phi_j(y)}r_j(y)$
for $j=1,\ldots,n$ and\/ $y\in\R$ by
\begin{gather}
r_j(y)=\sqrt{\ts\frac{1}{a_j}+y^2} \quad\text{and}\quad
\phi_j(y)=\psi_j+\int_0^y \frac{\d
t}{(\frac{1}{a_j}+t^2)\sqrt{P(t)}}\,,
\label{lm5eq1}\\
\text{where}\quad \ts
P(t)=\frac{1}{t^2}\bigl(\prod_{k=1}^n(1+a_kt^2)e^{\al t^2}-1\bigr).
\quad\text{Then}
\label{lm5eq2}\\
L=\bigl\{\bigl(x_1w_1(y),\ldots,x_n w_n(y)\bigr):
\text{$x_1,\ldots,x_n\in\R,$ $\ts\sum_{j=1}^nx_j^2=1$}\bigr\}
\label{lm5eq3}
\end{gather}
is a closed, embedded Lagrangian diffeomorphic to ${\mathcal
S}^{n-1}\t\R$ and satisfying\/ $\al F^\perp=H$. If\/ $\al>0$ it is a
self-expander, and if\/ $\al=0$ it is one of Lawlor's examples of
special Lagrangian submanifolds {\rm\cite{Lawl}}. It has Lagrangian
angle}
\begin{equation}
\th(y)=\ts\sum_{j=1}^n\phi_j(y)+\arg\bigl(y+iP(y)^{-1/2}\bigr).
\label{lm5eq4}
\end{equation}

\begin{proof} Suppose we are in the situation of Theorems A and B,
with $\la_1=\cdots=\la_n=C=1$, $\al\ge 0$, $A<0$ and $I=I_{\max}$.
Define $G(u)=Q(u)e^{\al u}$ and $\be=-\min_{j=1,\ldots,n}
\frac{\al_j}{\la_j}<0$. Then $G(\be)=0$, and
\begin{equation}
\frac{\d}{\d u}\ln(G(u))=\sum_{j=1}^n\frac{\la_j}{\al_j+\la_ju}+\al,
\label{lm5eq5}
\end{equation}
which is positive for $u>\be$. Therefore $G$ is an increasing
function on $[\be,\iy)$ with $G(\be)=0$ and $\lim_{u\ra
\iy}G(u)=\iy$. Note that $G(0)=\prod_{j=1}^n\al_j\ge A^2$, since
\eq{lm4eq9} at $s=s_0$ gives $(\prod_{j=1}^n\al_j)^{1/2}
\sin(\phi-\th)=A$ and $\md{\sin(\phi-\th)}\le 1$. Hence there exists
$u_*\in(\be,0]$ with $G(u_*)=A^2$ (see Figure~\ref{lm5fig}).

\begin{figure}[htb]
\begin{center}
\resizebox{5cm}{!}{\includegraphics{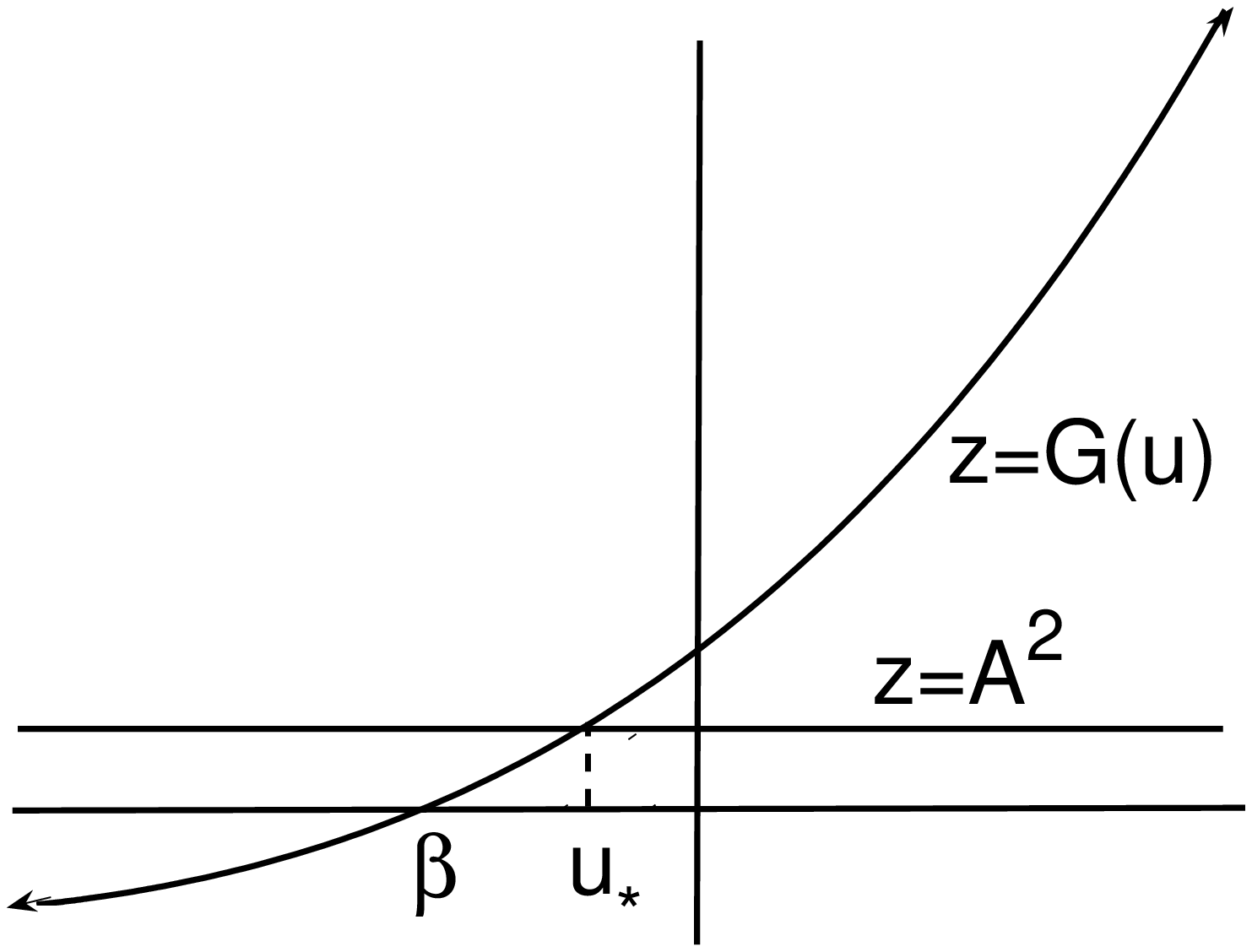}}
\end{center}
\vskip -1cm \caption{the case $\la_1=\cdots=\la_n=C=1$, $\al\ge
0$, and $A<0$.}
\label{lm5fig}
\end{figure}

Since $\md{\sin(\phi-\th)}\le 1$, equation \eq{lm4eq9} implies
that $G(u(s))\ge A^2$ for all $s\in I$, and so $u(s)\ge u_*$ for
all $s\in I$. Also, combining \eq{lm4eq7} and \eq{lm4eq9}, we have
$\frac{\d u}{\d s}=0$ $\Leftrightarrow$ $\cos(\phi-\th)=0$
$\Leftrightarrow$ $\md{\sin(\phi-\th)}=1$ $\Leftrightarrow$
$G(u(s))=A^2$ $\Leftrightarrow$ $u(s)=u_*$. If $u$ is a constant
function,  then $u(s)\equiv 0$ as we have $u(s_0)=0$, and also
$\phi-\th \equiv -\frac{\pi}{2}$ as $A<0$. From \eq{lm4eq7} and
\eq{lm4eq9}, it follows that $ \phi_j =\psi_j-\frac{As}{\al_j}$
and $\th=\th_0+\al As$ for some constants $\psi_j$ and $\th_0$. As
$\al_j>0$ and $\al \ge 0$, it contradicts to the fact that
$\phi-\th \equiv -\frac{\pi}{2}$. Now suppose $s_1,s_2$ are
distinct, adjacent zeroes of $\frac{\d u}{\d s}$ in $I$. Then
$\frac{\d u}{\d s}$ has constant sign in $(s_1,s_2)$, but
$u(s_1)=u(s_2)=u_*$, giving a contradiction by the Mean Value
Theorem. Hence $\frac{\d u}{\d s}$ has at most one zero in~$I$.

Write $I=(a,b)$ for $a,b\in\R\cup\{\pm\iy\}$. We claim that
$u(s)\ra\iy$ as $s\ra a_+$ or $s\ra b_-$. When $a$ or $b$ are finite
this follows from $I=I_{\max}$, since the only way the solution
could not extend over $a$ or $b$ is if $u\ra\iy$. When $a$ or $b$
are infinite, $u$ is monotone near infinity as $\frac{\d u}{\d s}$
has at most one zero, so $u(s)$ must approach infinity or some
finite limiting value $u'$ as $s\ra\pm\iy$. If $u(s)\ra u'$ as
$s\ra\pm\iy$ then $\frac{\d u}{\d s}\ra 0$ as $s\ra\iy$, forcing
$u'=u_*$ from above. We can exclude this possibility by showing that
$\frac{\d^2u}{\d s^2}\not\ra 0$ as $s\ra\iy$. Thus $u(s)\ra\iy$ as
$s\ra a_+$ or $s\ra b_-$, and $u$ has at least one minimum at $s_*$
in $I$. Then $\frac{\d u}{\d s}(s_*)=0$, so $s_*$ is unique from
above, and~$u(s_*)=u_*$.

Combining the first equation of \eq{lm4eq7}, \eq{lm4eq9}, and the
fact that $u$ is decreasing in $(a,s_*)$ and increasing in $(s_*,b)$
we have
\begin{equation}
\frac{\d u}{\d s}=\begin{cases}
-\,2\sqrt{Q(u)-A^2e^{-\al u}}\,, & a<s\le s_*,\\
\phantom{-}\,2\sqrt{Q(u)-A^2e^{-\al u}}\,, & s_*\le s<b.
\end{cases}
\label{lm5eq6}
\end{equation}
This gives
\begin{equation}
\int_{u_*}^{u(s)}\frac{\d v}{2\sqrt{Q(v)-A^2e^{-\al v}}}=\md{s-s_*}.
\label{lm5eq7}
\end{equation}
Thus $a,b$ are finite when $n > 2$ and $a=-\iy$, $b=+\iy$
when~$n\le 2$.

Equation \eq{lm5eq7} defines $s$ explicitly as a function of $u$.
Inverting this gives $u$ as a function of $s$. Then from \eq{lm4eq9}
we obtain $\sin(\phi-\th)$, and hence $\cos(\phi-\th)$, as functions
of $s$. Thus we have the right hand side of each equation in
\eq{lm4eq7} as functions of $s$, and integrating \eq{lm4eq7} gives
$\phi_j,\phi$ and $\th$ as functions of~$s$.

To make this more explicit, not depending on inverting the
integral function~\eq{lm5eq7}, we shall change from $s$ to a new
variable $y$ defined by
\begin{equation*}
\begin{split}
y(s) &=
\begin{cases}
-\sqrt{u-u(s_*)}, & a<s\le s_*,\\
\phantom{-}\sqrt{u-u(s_*)}, & s_*\le s<b.
\end{cases}
\end{split}
\end{equation*}
Then $y:(a,b)\ra\R$ is a smooth diffeomorphism, and $u=u(s_*)+y^2$,
so that $r_j^2=\al_j+u=\al_j+u(s_*)+y^2=a_j^{-1}+y^2$, where
$a_j=(\al_j+u(s_*))^{-1}$. This gives
$r_j(y)=\sqrt{\frac{1}{a_j}+y^2}$, as in \eq{lm5eq1}. Computing
$\frac{\d\phi_j}{\d y}$ from $\frac{\d\phi_j}{\d s}$ and $\frac{\d
y}{\d s}$ yields $\frac{\d\phi_j}{\d y}=\bigl((\frac{1}{a_j}+y^2)
\sqrt{P(y)}\bigr)^{-1}$, for $P(y)$ as in \eq{lm5eq2}. This implies
the second equation of \eq{lm5eq1}, with $\psi_j=\phi_j\vert_{y=0}
=\phi_j\vert_{s=s_*}$.

Theorems A and B now imply that $L$ is Lagrangian with $\al
F^\perp=H$, with Lagrangian angle \eq{lm5eq4}. Equation
\eq{lm5eq3} implies that $L$ is diffeomorphic to ${\mathcal
S}^{n-1}\t\R$, and closedness of $L$ follows from
$r_j(y)\ra\pm\iy$ as $y\ra\pm\iy$. That $L$ is embedded follows
from the fact that each $\phi_j(y)$ is strictly increasing, and
has image an interval of size at most $\pi$, as we will show in
the proof of Theorem D. When $\al=0$ our formulae reduce to those
of Harvey's treatment \cite[p.~139--143]{Harv} of Lawlor's
examples \cite{Lawl}. This completes the proof.
\end{proof}

{\bf Theorem D.} {\it In the situation of Theorem C, there exist\/
$\bar\phi_1,\ldots,\bar\phi_n\in(0,\frac{\pi}{2}]$ with\/
$\bar\phi_j=\int_0^\infty \frac{\d
t}{(\frac{1}{a_j}+t^2)\sqrt{P(t)}}$\/ for\/ $ j=1,\ldots,n,$ such
that the Lagrangian $L$ is asymptotic at infinity to the union of
Lagrangian planes $L_1 \cup L_2,$ where
\begin{equation}
\begin{split}
L_1&=\bigl\{(e^{i(\psi_1+\bar{\phi}_1)}t_1,\ldots,
e^{i(\psi_n+\bar{\phi}_n)}t_n):t_1,\ldots,t_n\in\R\bigr\},\\
L_2&=\bigl\{(e^{i(\psi_1-\bar{\phi}_1)}t_1,\ldots,
e^{i(\psi_n-\bar{\phi}_n)}t_n):t_1,\ldots,t_n\in\R\bigr\}.
\end{split}
\label{lm5eq8}
\end{equation}
We have $0<\bar\phi_1+\cdots+\bar\phi_n<\frac{\pi}{2}$ if\/ $\al>0,$
and\/ $\bar\phi_1+\cdots+\bar\phi_n=\frac{\pi}{2}$ if\/~$\al=0$.

Fix $\al>0$. Then $\Phi^n:(a_1,\ldots,a_n)\mapsto(\bar\phi_1,
\ldots,\bar\phi_n)$ gives a diffeomorphism
\begin{equation}
\Phi^n:(0,\iy)^n\longra\bigl\{(\bar\phi_1,\ldots,\bar\phi_n)\in(0,
\ts\frac{\pi}{2}]^n:
0<\bar\phi_1+\cdots+\bar\phi_n<\frac{\pi}{2}\bigr\}.
\label{lm5eq9}
\end{equation}
That is, for all\/ $\al>0$ and\/ $L_1,L_2$ satisfying
$0<\bar\phi_1+\cdots+\bar\phi_n<\frac{\pi}{2}$ as above, Theorem C
gives a unique Lagrangian expander $L$ asymptotic to~$L_1\cup L_2$.

When $\al=0,$ it is studied by  Lawlor in {\rm \cite{Lawl}}. The map
$\Phi^n:(a_1,\ldots,a_n)\ab\mapsto(\bar\phi_1, \ldots,\bar\phi_n)$
gives a surjection
\begin{equation}
\Phi^n:(0,\iy)^n\longra\bigl\{(\bar\phi_1,\ldots,\bar\phi_n)\in(0,
\ts\frac{\pi}{2})^n:
\bar\phi_1+\cdots+\bar\phi_n=\frac{\pi}{2}\bigr\}, \label{lm5eq10}
\end{equation}
such that\/ $(a_1,\ldots,a_n)$ and\/ $(a_1',\ldots,a_n')$ have the
same image $(\bar\phi_1,\ldots,\bar\phi_n)$ if and only if\/
$a_j'=ta_j$ for some $t>0$ and all\/ $j=1,\ldots,n,$ and the
corresponding special Lagrangians $L,L'$ satisfy~$L'=t^{-1/2}L$. }
\medskip

\begin{proof} From the definition of $\phi_{j}(y)$ in \eq{lm5eq1},
it is clear that the integral converges as $y\ra\iy$ which is
denoted by $\bar{\phi}_j>0$. Here $\bar{\phi}_j$ depends on
 $a_1,\ldots,a_n>0$ and $\al\ge 0$. The limit of the integral as
$y\ra-\iy$ is then $-\bar{\phi}_j$. This shows that $L$ is
asymptotic to $L_1\cup L_2$. It is also easy to see that when
$\al>0$, the convergence of $L$ to $L_1\cup L_2$ at infinity is
exponential.

Since $P(y)^{-1/2}>0$ and $\md{y}\gg P(y)^{-1/2}$ for large $y$ by
\eq{lm5eq2}, we see that $\lim_{y\ra-\iy}\arg\bigl(y+iP(y)^{-1/2}
\bigr)=\pi$ and $\lim_{y\ra\iy}\arg\bigl(y+iP(y)^{-1/2}\bigr)=0$.
Thus \eq{lm5eq4} implies that
\begin{equation}
\ts
\lim\limits_{y\ra-\iy}\th(y)=\sum\limits_{j=1}^n\psi_j
-\sum\limits_{j=1}^n\bar\phi_j+\pi,\;\>
\lim\limits_{y\ra\iy}\th(y)=\sum\limits_{j=1}^n\psi_j
+\sum\limits_{j=1}^n\bar\phi_j.
\label{lm5eq11}
\end{equation}
But the last equation of \eq{lm4eq7}, \eq{lm4eq9}, and $A<0$ imply
that $\th$ is strictly decreasing when $\al>0$ and constant when
$\al=0$. Hence $\lim_{y\ra-\iy}\th(y)>\lim_{y\ra\iy}\th(y)$ when
$\al>0$ and $\lim_{y\ra-\iy}\th(y)=\lim_{y\ra\iy}\th(y)$ when
$\al=0$. By \eq{lm5eq11}, this gives
$\sum_{j=1}^n\bar\phi_j<\frac{\pi}{2}$ when $\al>0$, and
$\sum_{j=1}^n\bar\phi_j=\frac{\pi}{2}$ when $\al=0$. As each
$\bar\phi_k>0$, this implies that $\bar\phi_j\le\frac{\pi}{2}$,
and completes the first part.

Write the map $\Phi^n$ in the theorem as $\Phi^n=(\Phi^n_1,\ldots,
\Phi^n_n)$. Then \eq{lm5eq1}--\eq{lm5eq2} and the definition of
$\Phi^n$ imply that
\begin{equation}
\Phi^n_j(a_1,\ldots,a_n)=\int_0^\iy\frac{a_j\,\d
y}{(1+a_jy^2)\sqrt{
\frac{1}{y^2}\bigl(\prod_{l=1}^n(1+a_ly^2)e^{\al y^2}-1\bigr)}}\,.
\label{lm5eq12}
\end{equation}
Computation shows that for $y\in(0,\iy)$ we have
\begin{equation*}
\frac{\pd}{\pd a_k}
\raisebox{-3pt}{\begin{Large}$\displaystyle\Bigl[$\end{Large}}
\frac{a_j}{(1+a_jy^2)\sqrt{
\frac{1}{y^2}\bigl(\prod_{l=1}^n(1+a_ly^2)e^{\al y^2}-1\bigr)}}
\raisebox{-3pt}{\begin{Large}$\displaystyle\Bigr]$\end{Large}}
<0\quad\text{for $k\ne j$.}
\end{equation*}
Integrating this over $(0,\iy)$ and using \eq{lm5eq12} thus gives
\begin{equation}
\ts\frac{\pd}{\pd
a_k}\bigl(\Phi^n_j(a_1,\ldots,a_n)\bigr)<0\;\>\text{for $k\ne j$.}
\label{lm5eq13}
\end{equation}
Let $t>0$. Replacing $a_j$ by $ta_j$ for $j=1,\ldots,m$ in
\eq{lm5eq12} and changing variables from $y$ to $t^{-1/2}y$ shows
that
\begin{equation}
\Phi^n_j(ta_1,\ldots,ta_n)=\int_0^\iy\frac{a_j\,\d
y}{(1+a_jy^2)\sqrt{
\frac{1}{y^2}\bigl(\prod_{l=1}^n(1+a_ly^2)e^{t^{-1}\al
y^2}-1\bigr)}}\,. \label{lm5eq14}
\end{equation}
The integrand here is a strictly increasing function of $t$ for
$\al>0$, and constant for $\al=0$. Thus, $\frac{\d}{\d
t}\Phi^n_j(ta_1,\ldots, ta_n)$ is positive for $\al>0$ and zero for
$\al=0$. Setting $t=1$ yields
\begin{equation}
\sum_{k=1}^{n}a_k\frac{\pd}{\pd
a_k}\bigl(\Phi^n_j(a_1,\ldots,a_n)\bigr)
\begin{cases} >0, & \al>0, \\ =0, & \al=0. \end{cases}
\label{lm5eq15}
\end{equation}
Combining \eq{lm5eq13} and \eq{lm5eq15} implies that
\begin{equation}
\frac{\pd}{\pd a_j}\bigl(\Phi^n_j(a_1,\ldots,a_n)\bigr)>0.
\label{lm5eq16}
\end{equation}
Also, taking the limit $t\ra\iy$ in \eq{lm5eq14} we see that
$\lim_{t\ra\iy}\Phi^n_j(ta_1,\ldots,ta_n)$ exists, and equals
$\bar\phi_j$ with the same $a_1,\ldots,a_n$ but with $\al=0$. But we
have already shown that $\bar\phi_1+\cdots+\bar\phi_n=\frac{\pi}{2}$
when $\al=0$. Therefore
\begin{equation}
\ts\sum_{j=1}^n\lim_{t\ra\iy}\Phi^n_j(ta_1,\ldots,ta_n)=\frac{\pi}{2}.
\label{lm5eq17}
\end{equation}

Fixing $\al>0$, we first show that the differential of $\Phi^n$ is
nonsingular. Suppose there exist $\la_1,\ldots,\la_n\in\R$ not all
zero such that for $j=1,\ldots,n$ we have
$\sum_{k=1}^n\la_k\frac{\pd}{\pd
a_k}\bigl(\Phi^n_j(a_1,\ldots,a_n) \bigr)=0$. Pick $j$ such that
$\md{\la_j}/a_j$ is greatest. Then \eq{lm5eq13} and \eq{lm5eq15}
imply that
\begin{equation*}
\begin{split}
{\ts a_j\frac{\pd}{\pd a_j}\bigl(\Phi^n_j(a_1,\ldots,a_n)\bigr)}&>
{\ts-\sum_{k\ne j}a_k\frac{\pd}{\pd
a_k}\bigl(\Phi^n_j(a_1,\ldots,a_n)\bigr)} \quad\text{and}\\
{\ts \frac{|\la_j|}{a_j}\,a_j\frac{\pd}{\pd
a_j}\bigl(\Phi^n_j(a_1,\ldots,a_n)\bigr)}&> {\ts-\sum_{k\ne
j}\frac{|\la_k|}{a_k}\,a_k\frac{\pd}{\pd
a_k}\bigl(\Phi^n_j(a_1,\ldots,a_n)\bigr)}.
\end{split}
\end{equation*}
It follows that
\begin{equation*}
\bmd{\ts\la_j\frac{\pd}{\pd
a_j}\bigl(\Phi^n_j(a_1,\ldots,a_n)\bigr)}> \bmd{\ts\sum_{k\ne
j}\la_k\frac{\pd}{\pd a_k}\bigl(\Phi^n_j(a_1,\ldots,a_n)\bigr)},
\end{equation*}
contradicting $\sum_{k=1}^n\la_k\frac{\pd}{\pd
a_k}\bigl(\Phi^n_j(a_1,\ldots,a_n)\bigr)=0$. Thus no such
$\la_1,\ldots,\la_n$ exist, and
$\d\Phi^n\vert_{(a_1,\ldots,a_n)}:\R^n\ra\R^n$ is invertible. So
$\Phi^n$ in \eq{lm5eq9} is a local diffeomorphism. The same argument
when $\al=0$ shows that the only possible $(\la_1,\ldots,\la_n)$ are
multiples of $(a_1,\ldots,a_n)$. So in \eq{lm5eq10},
$\d\Phi^n\vert_{(a_1,\ldots,a_n)}:\R^n\ra\R^{n-1}$ has kernel
$\an{(a_1,\ldots,a_n)}$, and is surjective.

We will now show that when $\al>0$, the map $\Phi^n$ of
\eq{lm5eq9} is surjective. Embed the domain $(0,\iy)^n$ of
$\Phi^n$ in $\mathbb{RP}^n$ by $(a_1,\ldots,a_n)\mapsto
[1,a_1,\ldots,a_n]$. The closure of $(0,\iy)^n$ in $\mathbb{RP}^n$
is an $n$-simplex $\De^n$. It has faces $\De^{n-1}_j$ for
$j=0,\ldots,n$ on which $x_j=0$ in homogeneous coordinates
$[x_0,\ldots,x_n]$. The closure in $\R^n$ of the range of $\Phi^n$
in \eq{lm5eq9} is also an $n$-simplex $\ti\De^n$ with faces
$\ti\De^{n-1}_j$ for $j=0,\ldots,n$, where $\bar\phi_1+\cdots+
\bar\phi_n=\frac{\pi}{2}$ on~$\ti\De^{n-1}_0$ and $\bar\phi_j=0$
on $\ti\De^{n-1}_j$ for $j=1,\ldots,n$. Note that
$\pd\De^n=\cup_{j=0}^n \De^{n-1}_j$ and $\pd\ti\De^n=\cup_{j=0}^n
\ti\De^{n-1}_j$.

We claim $\Phi^n$ extends to a continuous map
$\bar\Phi^n:\De^n\ra\ti\De^n$, which maps
$\De^{n-1}_j\ra\ti\De^{n-1}_j$ for $j=0,\ldots,n$. To see this, note
that $\lim_{a_j\ra 0}\Phi^n_k(a_1,\ab\ldots,\ab a_n)$ exists for all
$k$ by \eq{lm5eq12}, and $\lim_{a_j\ra
0}\Phi^n_j(a_1,\ldots,a_n)=0$. Thus $\Phi^n$ extends continuously to
$\De^{n-1}_j$ for $j=1,\ldots,n$, and maps
$\De^{n-1}_j\ra\ti\De^{n-1}_j$. Also, the fact above that
$\lim_{t\ra\iy}\Phi^n_j(ta_1,\ldots,ta_n)$ exists shows that
$\Phi^n$ extends to $\De^{n-1}_0$, and \eq{lm5eq17} implies that
this extension maps~$\De^{n-1}_0\ra\ti\De^{n-1}_0$.

We will prove surjectivity of \eq{lm5eq9} and its extension
$\bar\Phi^n$ by induction on $n$. The map is clearly onto when
$n=1$ since it is continuous and $\bar\Phi^1([1,0])=0$,
$\bar\Phi^1([0,1])=\frac{\pi}{2}$ by \eq{lm5eq17}. Suppose
$\bar\Phi^{n-1}$ is surjective. Since $\bar\Phi^n$ reduces to
$\bar\Phi^{n-1}$ when $a_k=0$, this implies that
$\bar\Phi^n\vert_{\De^{n-1}_k}:\De^{n-1}_k\ra\ti\De^{n-1}_k$ is
surjective for $k=1,\ldots,n$. Now consider
$\bar\Phi^n\vert_{\De^{n-1}_0}:\De^{n-1}_0 \ra\ti\De^{n-1}_0$.
Since $\bar\Phi^n\vert_{\De^{n-1}_k}$ is surjective for
$k=1,\ldots,n$, we see that
$\bar\Phi^n\vert_{\De^{n-1}_0\cap\De^{n-1}_k}:
\De^{n-1}_0\cap\De^{n-1}_k\ra\ti\De^{n-1}_0\cap\ti\De^{n-1}_k$ is
surjective for $k=1,\ldots,n$. So $\bar\Phi^n$ takes
$\pd\De^{n-1}_0$ surjectively to $\pd\ti\De^{n-1}_0$, and is of
degree one. Using algebraic topology, it follows that
$\bar\Phi^n\vert_{\smash{\De^{n-1}_0}}:
\De^{n-1}_0\ra\ti\De^{n-1}_0$ is surjective. Hence $\bar\Phi^n$
takes $\pd\De^n$ surjectively to $\pd\ti\De^n$, and is of degree
one, so again, $\Phi^n$ is surjective.

Therefore by induction, $\Phi^n$ in \eq{lm5eq9} is surjective for
all $n$. But $\Phi^n$ is a local diffeomorphism, and extends to a
map $\De^n\ra\ti\De^n$ taking $\pd\De^n\ra\pd\ti\De^n$, so $\Phi^n$
is proper, and thus $\Phi^n$ is a covering map. As the domain of
$\Phi^n$ is connected and the range simply-connected, it follows
that $\Phi^n$ in \eq{lm5eq9} is a {\it diffeomorphism}, as we have
to prove. The final part for \eq{lm5eq10} follows by a similar
argument; one way to do it is to show that the restriction of
$\Phi^n$ in \eq{lm5eq10} to $\bigl\{(a_1,\ldots,a_n)
\in(0,\iy)^n:a_1+\cdots+a_n=1\bigr\}$ is a diffeomorphism.
\end{proof}

In the last two parts of Theorem D, the proof that $\Phi^n$ is
surjective is based on Lawlor~\cite[Lemma 10]{Lawl}.

\section{Other self-similar solutions}
\label{lm6}

Finally we prove Theorems E and F of~\S\ref{lm33}.
\medskip

{\bf Theorem E.} {\it In Theorems A and B, suppose that either:
\begin{itemize}
\setlength{\itemsep}{0pt}
\setlength{\parsep}{0pt}
\item[{\rm(a)}]  $\la_1=\cdots=\la_n=C=1,$ $\al<0$ and\/ $A>0;$ or
\item[{\rm(b)}] $\la_1=\cdots=\la_m=1$ and\/
$\la_{m+1}=\cdots=\la_n=-1$ for some $1\le m<n,$ $C=1,$ $A>0,$
and\/~$\al\in\R$.
\end{itemize}
Then solutions exist for all\/ $s\in\R,$ and we take\/ $I=\R$. In
each of cases {\rm(a),(b)} we divide into two subcases:
\begin{itemize}
\setlength{\itemsep}{0pt}
\setlength{\parsep}{0pt}
\item[{\rm(i)}] $\sum_{j=1}^n\frac{\la_j}{\al_j}+\al=0$ and\/
$\al_1\cdots\al_n=A^2;$ or
\item[{\rm(ii)}] otherwise.
\end{itemize}
In case {\rm(i),} we have explicit solutions to \eq{lm4eq7} and
obtain
\begin{equation}
\begin{split}
L=\bigl\{\bigl(x_1\sqrt{\al_1}\,&e^{i(\psi_1-\la_1As/\al_1)},
\ldots,x_n\sqrt{\al_n}\,e^{i(\psi_n-\la_nAs/\al_n)}\bigr):\\
&\text{$x_1,\ldots,x_n\in\R,$ $s\in\R,$
$\ts\sum_{j=1}^n\la_jx_j^2=1$}\bigr\},
\end{split}
\label{lm6eq2}
\end{equation}
which is Hamiltonian stationary in addition to being self-similar,
and invariant under a subgroup $\R$ or $\U(1)$ of diagonal matrices
$\bigl\{\diag(e^{i\la_1t/\al_1},\ab \ldots,\ab
e^{i\la_nt/\al_n}):t\in\R\bigr\}$ in~$\U(n)$.

In case {\rm(ii),} $u$ and\/ $\phi-\th$ are periodic in $s$ with
period\/ $S>0$, and
\begin{equation}
\begin{aligned}
u(s+S)&=u(s),& \phi_j(s+S)&=\phi_j(s)+\ga_j,\\
\phi(s+S)&=\phi(s)+\ts\sum_{j=1}^n\ga_j, &
\th(s+S)&=\th(s)+\ts\sum_{j=1}^n\ga_j,
\end{aligned}
\label{lm6eq3}
\end{equation}
for some $\ga_1,\ldots,\ga_n\in\R$ and all\/ $s\in\R$. In case
{\rm(b)} with\/ $\al=0$ we have $\th(s)\equiv\th(0)$
and\/~$\sum_{j=1}^n\ga_j=0$.}
\medskip

Before the proof of Theorem E, we first derive the following lemma
and proposition:

\begin{lem} In the situation of Theorem E, $G(u)=Q(u)e^{\al u}$
has a unique critical point $u^*$ on the interval $(\be_1,\be_2),$
where $\be_1<0<\be_2$ are defined by
\begin{equation}
\be_1 \!=\!\begin{cases}
-\min_{1\leq j\leq n}\al_j, & \!\!\!\text{in {\rm(a),}}\\
-\min_{1\leq j\leq m}\al_j, & \!\!\!\text{in {\rm(b),}}
\end{cases}\;\>
\be_2\!=\!\begin{cases}
\iy, & \!\!\!\text{in {\rm(a),}}\\
\min_{m+1\leq j\leq n}\al_j, & \!\!\!\text{in {\rm(b).}}
\end{cases}
\label{lm6eq4}
\end{equation}
Also $\lim_{u\ra\be_1}G(u)=\lim_{u\ra\be_2}G(u)=0$, $G(u)>0$ on
$(\be_1,\be_2)$, $G'(u)>0$ on $(\be_1,u^*)$, $G'(u)<0$ on
$(u^*,\be_2)$, and\/ $u(s)\in(\be_1,\be_2)$ for all\/~$s\in I$.
\label{lm6lem1}
\end{lem}

\begin{proof} The first derivative of $\ln G$ is given in
\eq{lm5eq5}. Differentiating yields
\begin{equation*}
\ts\frac{\d^2}{\d u^2}\ln(G(u))=-\sum_{j=1}^n\frac{\la_j^2}{(\al_j+
\la_ju)^2}<0.
\end{equation*}
Hence $\frac{\d}{\d u}\ln(G(u))$ is strictly decreasing, and
$\frac{\d}{\d u}\ln(G(u))$ can have at most one zero in any
interval on which $G(u)>0$ so $\ln(G(u))$ is defined. By
definition of $G$ we see that
$\lim_{u\ra\be_1}G(u)=\lim_{u\ra\be_2}G(u)=0$ and $G(u)>0$ on
$(\be_1,\be_2)$. Thus $G$ must have a global maximum $u^*$ in
$(\be_1,\be_2)$. Then $u^*$ is a zero of $\frac{\d}{\d
u}\ln(G(u))$ in $(\be_1,\be_2)$, so $u^*$ is unique, and $G'(u)>0$
on $(\be_1,u^*)$ and $G'(u)<0$ on $(u^*,\be_2)$ follow as
$\frac{\d^2}{\d u^2}\ln(G(u))<0$. Finally, since
$r_j^2(s)=\al_j+\la_ju(s)>0$ for all $s\in I$ and $j=1,\ldots,n$,
we see that $u(s)\in(\be_1,\be_2)$ for all\/ $s\in I$
from~\eq{lm6eq4}.
\end{proof}

\begin{prop} In case {\rm(ii),} there exist unique, finite
$u_1,u_2$ with\/ $\be_1<u_1<u^*<u_2<\be_2$ and\/ $G(u_1)=G(u_2)=A^2$
and\/ $G(u)>A^2$ on $(u_1,u_2)$. We have $u(s)\in[u_1,u_2]$ for
all\/ $s\in I,$ and solutions exist for all\/~$s\in\R$.
\label{lm6prop1}
\end{prop}

\begin{proof} Since $\sin^2(\phi-\th)\le 1$, equation \eq{lm4eq9}
implies that $G(u(s))\ge A^2$ for all $s\in I$. By Lemma
\ref{lm6lem1}, if $s\in I$ then $u(s)\in(\be_1,\be_2)$, so
$G(u(s))\le G(u^*)$. Thus $G(u^*)\ge A^2$. If $G(u^*)=A^2$ this
forces $u(s)=u^*$ for all $s\in I$, so $u^*=u(s_0)=0$, giving
$G(u^*)=\al_1\cdots\al_n=A^2$, and $\frac{\d}{\d
u}\ln(G(u))\vert_{u=0}=0$, giving $\sum_{j=1}^n\frac{\la_j}{\al_j}
+\al=0$ by \eq{lm5eq5}. Thus $G(u^*)=A^2$ implies we are in
case~(i).
\begin{figure}[htb]
\vskip -.2cm
\begin{center}
\resizebox{5cm}{!}{\includegraphics{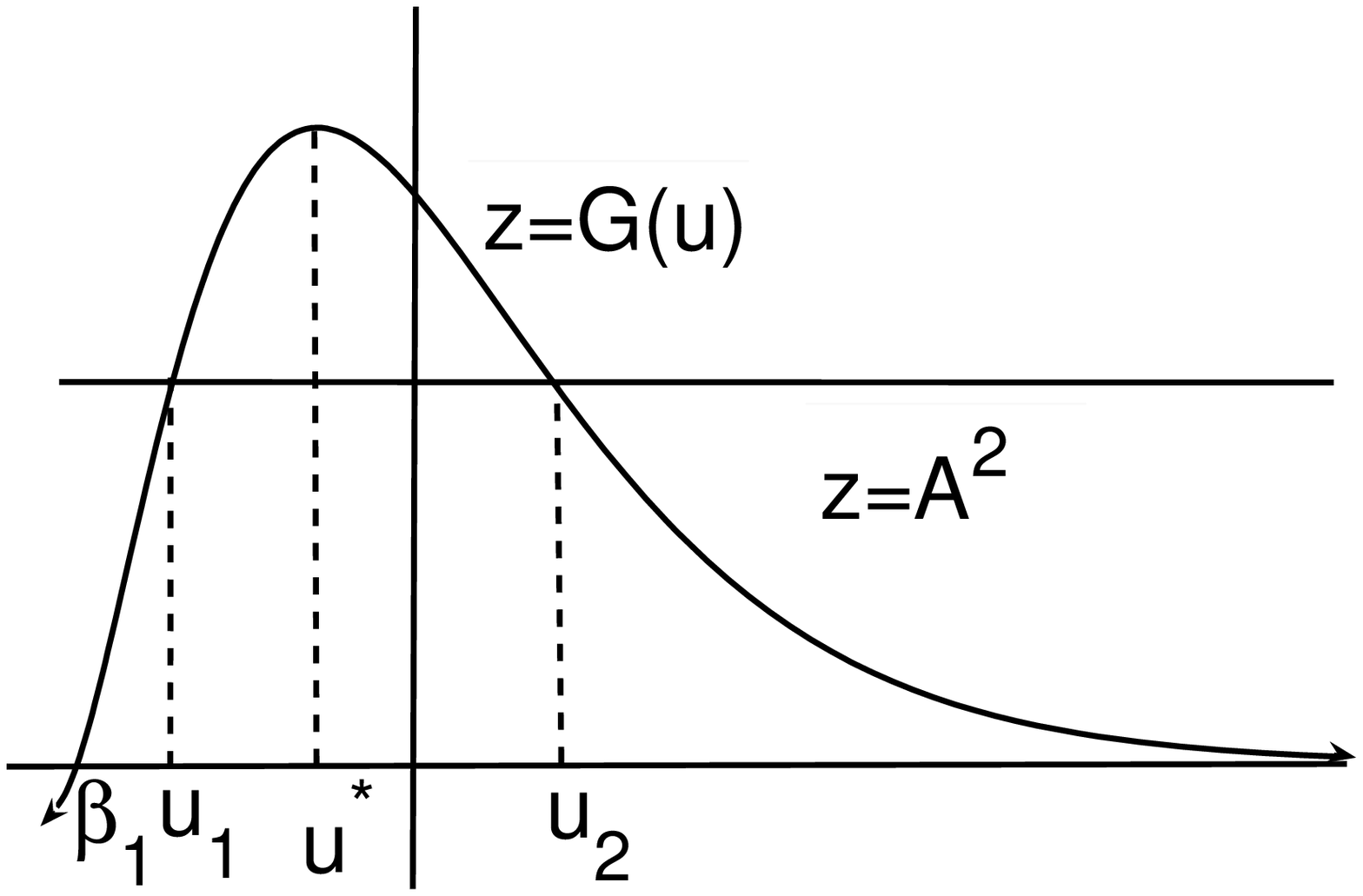}}
\resizebox{5cm}{!}{\includegraphics{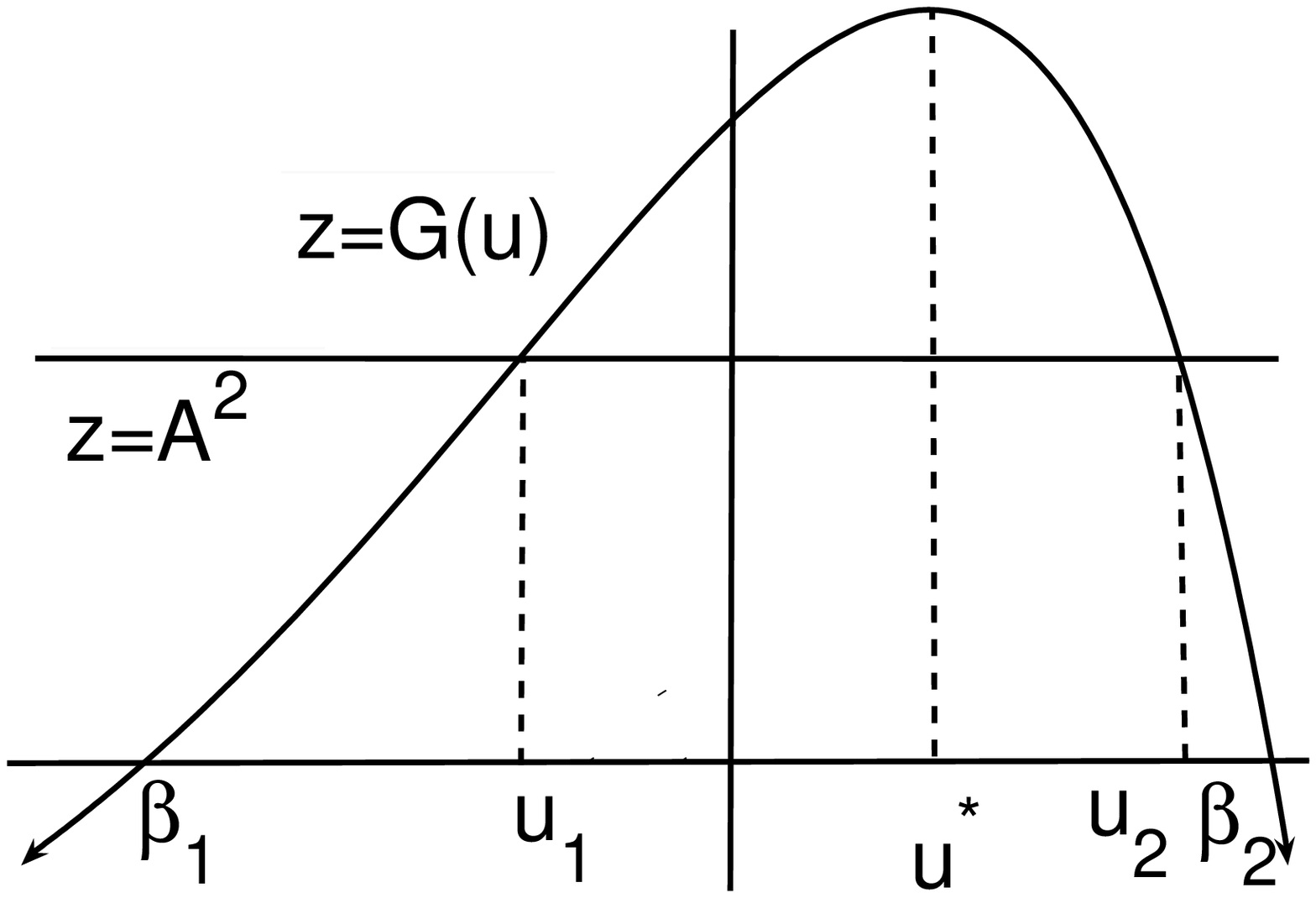}}
\end{center}
\vskip -.8cm \caption{Case (a)\hskip 3.5cm Case (b)}
\label{lm6fig}
\end{figure}

Since we restricted to case (ii), we have $G(u^*)>A^2$. So Lemma
\ref{lm6lem1} and the Intermediate Value Theorem imply that there
exist unique $u_1\in(\be_1,u^*)$ and $u_2\in(u^*,\be_2)$ with
$G(u_1)=G(u_2)=A^2$, and that if $u\in(\be_1,\be_2)$ then $G(u)\ge
A^2$ if and only if $u\in[u_1,u_2]$, with $G(u)>A^2$ on $(u_1,u_2)$.
Thus $G(u(s))\in[u_1,u_2]$ for all $s\in I$, by Lemma \ref{lm6lem1}.
This is illustrated in Figure~\ref{lm6fig}.

Suppose that in Theorem E, we have solutions on some interval $I$ in
$\R$. These must extend to some maximal open interval
$I_{\max}=(a,b)$ in $\R$, for $a,b\in\R\cup\{\pm\iy\}$. We could
only have $a>-\iy$ if either $u(s)\ra\iy$ as $s\ra a_+$ or
$\al_j+\la_ju(s)\ra 0$ as $s\ra a_+$ for some $j=1,\ldots,n$, so
that the right hand side of some equation in \eq{lm4eq7} becomes
singular as $s\ra a_+$, and the solutions do not extend past $a$.
But this is impossible because it can only happen when $u(s)$
approaches $\be_1$ or $\be_2$, and  $u(s)$ is confined to
$[u_1,u_2]$ which lies strictly inside $(\be_1,\be_2)$ from the
discussion above. Therefore $a=-\iy$, and similarly $b=\iy$, and
solutions exist for all~$s\in\R$.
\end{proof}

Now we are ready to prove Theorem E.

\begin{proof}[Proof of Theorem E]
In case (i), $G(0)=G(u^*)=A^2$, equation \eq{lm4eq9} implies
$\sin^2(\phi-\th)\equiv 1$. It is easy to verify that solutions to
\eq{lm4eq7} are of the form
\begin{equation}
\begin{gathered}
u(s)\!=\!0,\;\, Q(u(s))\!=\!\al_1\cdots\al_n=A^2,\;\, \ts
\th(s)\!=\!\ts\sum_{j=1}^n\psi_j\!-\!\frac{\pi}{2}\!+\!\al A s,\\
\ts\phi_j(s)\!=\!\psi_j\!-\!\frac{\la_jAs}{\al_j},\;\,
\phi(s)\!=\!\sum_{j=1}^n\psi_j\!-\!\sum_{j=1}^n\frac{\la_jAs}{\al_j}
\!=\!\sum_{j=1}^n\psi_j\!+\!\al A s
\end{gathered}
\label{lm6eq1}
\end{equation}
for some $\psi_1,\ldots,\psi_n\in\R,$  \/ which exist for all
$s\in\R$. \/ From \eq{lm4eq8} we obtain $L$ as described in
\eq{lm6eq2} and have the induced metric on $L$ described in
\eq{lm4eq3}. Therefore, $\th$ is harmonic and $L$ is Hamiltonian
stationary.

In case (ii), we already proved that solutions exist for all
$s\in\R$ in Proposition \ref{lm6prop1}. It remains to show the
solutions are periodic. The proof of \eq{lm5eq6} implies that
\begin{equation*}
\ts\bigl(\frac{\d u}{\d s}\bigr)^2=4(Q(u)-A^2e^{-\al
u}\bigr)=4e^{-\al u}\bigl(G(u)-A^2\bigr).
\end{equation*}
Thus $\frac{\d u}{\d s}=0$ if and only if $G(u)=A^2$, that is, if
and only if $u=u_1$ or $u=u_2$. So $\frac{\d u}{\d s}$ cannot change
sign except at $s$ with $u(s)=u_1$ or $u(s)=u_2$, and $\frac{\d
u}{\d s}$ is determined up to sign by $u(s)$. As for \eq{lm5eq7},
the interval in $s$ taken for $u(s)$ to increase from $u_1$ to
$u_2$, or to decrease from $u_2$ to $u_1$, is
\begin{equation*}
\frac{S}{2}=\int_{u_1}^{u_2}\frac{\d v}{2\sqrt{Q(v)-A^2e^{-\al
v}}}\,,
\end{equation*}
which is finite, since $G'(u_1)>0$ and $G'(u_2)<0$ by Lemma
\ref{lm6lem1}, so $Q(v)-A^2e^{-\al v}$ has only simple zeroes at
$v=u_1$ and~$v=u_2$.

Therefore $u$ is periodic with period $S>0$, as it must increase
from $u=u_1$ to $u=u_2$ in an interval $S/2$, then decrease back to
$u=u_1$ in an interval $S/2$, and repeat. Hence $\frac{\d u}{\d s}$
has period $S$, so $\cos(\phi-\th)$ is periodic with period $S$ by
\eq{lm4eq7}. Thus $\phi-\th$ changes by an integral multiple of
$2\pi$ over each interval $S$. But \eq{lm4eq9} implies that
$\sin(\phi-\th)>0$, and $\phi-\th$ is continuous, so this multiple
of $2\pi$ is zero, and $\phi-\th$ is periodic with period~$S$.

Equation \eq{lm4eq7} now implies that $\frac{\d\phi_j}{\d s}$ is
periodic with period $S$. Integrating gives
$\phi_j(s+S)=\phi_j(s)+\ga_j$ for all $j=1,\ldots,n$ and $s\in\R$,
where $\ga_j=\int_0^S\frac{\d\phi_j}{\d s}(s)\d s$. Summing over
$j=1,\ldots,n$ gives $\phi(s+S)=\phi(s)+\sum_{j=1}^n\ga_j$ for all
$s\in\R$. Since $\phi-\th$ is periodic with period $S$ this
implies that $\th(s+S)=\th(s)+\sum_{j=1}^n\ga_j$, proving
\eq{lm6eq3}. In case (b) with $\al=0$ we have $\frac{\d\th}{\d
s}\equiv 0$ by \eq{lm4eq7}, so $\th(s)\equiv\th(0)$, and
$\sum_{j=1}^n\ga_j=0$. This completes the proof of Theorem~E.
\end{proof}

{\bf Theorem F.} {\it In Theorem E, we say that\/ $(w_1,\ldots,w_n)$
is periodic if there exists $T>0$ with\/ $w_j(s)=w_j(s+T)$ for all\/
$s\in\R$ and\/~$j=1,\ldots,n$.

If\/ $(w_1,\ldots,w_n)$ is periodic then in case {\rm(a),} $L$ is a
compact, immersed Lagrangian self-shrinker diffeomorphic to ${\mathcal
S}^1\t{\mathcal S}^{n-1},$ and in case {\rm(b),} $L$ is a closed,
noncompact, immersed Lagrangian diffeomorphic to ${\mathcal S}^1\t{\mathcal
S}^{m-1}\t\R^{n-m},$ a self-expander if\/ $\al>0,$ a self-shrinker
if\/ $\al<0,$ and special Lagrangian if\/~$\al=0$.

In case {\rm(i),} $(w_1,\ldots,w_n)$ is periodic if and only if\/
$\frac{\la_j}{\al_j}=\mu q_j$ with\/ $\mu>0$ and\/ $q_j\in\Q$ for
$j=1,\ldots,n$. In case {\rm(ii),} $(w_1,\ldots,w_n)$ is periodic if
and only if\/ $\ga_j\in\pi\Q$ for $j=1,\ldots,n$. In both cases, for
fixed\/ $m,\al,$ there is a dense subset of initial data for which\/
$(w_1,\ldots,w_n)$ is periodic.}
\medskip

\begin{proof} The first parts are straightforward. If
$(w_1,\ldots,w_n)$ is periodic with period $T$ then $L$ is the image
of an immersion $Q\t\R/T\Z\ra\C^n$, where $Q$ is the quadric
$\bigl\{(x_1,\ldots,x_n)\in\R^n:x_1^2+\cdots+x_n^2=1\bigr\}$ in (a),
which is diffeomorphic to ${\mathcal S}^{n-1}$, and the quadric
$\bigl\{(x_1,\ldots,x_n)\in\R^n:x_1^2+\cdots+x_m^2-x_{m+1}^2-\cdots
-x_n^2=1\bigr\}$ in (b), which is diffeomorphic to ${\mathcal
S}^{m-1}\t\R^{n-m}$. Since $\R/T\Z$ is diffeomorphic to ${\mathcal
S}^1$, $L$ is diffeomorphic as an immersed submanifold to ${\mathcal
S}^1\t{\mathcal S}^{n-1}$ in (a), which is compact, and to
${\mathcal S}^1\t{\mathcal S}^{m-1}\t\R^{n-m}$ in (b), which is
noncompact. It is a self-expander if $\al>0$, a self-shrinker if
$\al<0$, and special Lagrangian if~$\al=0$. We can also easily
verify that $L$ is closed when $(w_1,\ldots,w_n)$ is periodic.

The necessary and sufficient conditions for periodicity in the
last part are also easy. In case (i), if $(w_1,\ldots,w_n)$ is
periodic with period $T$ then $e^{i\phi_j(s+T)}=e^{i\phi_j(s)}$
for $j=1,\ldots,n$, so \eq{lm6eq1} gives $\frac{\la_jAT}{\al_j}\in
2\pi\Z$ for $j=1,\ldots,n$, and the condition holds with
$\mu=\frac{2\pi}{AT}>0$ and $q_j=\frac{\la_jAT}{2\pi\al_j}\in\Z
\subset\Q$. Conversely, if $\frac{\la_j}{\al_j}=\mu q_j$ for
$\mu>0$ and $q_j\in\Q$ then we may write $q_j=p_j/r$ for
$j=1,\ldots,n$, $p_j\in\Z$ and $r\in\N$ the lowest common
denominator of $q_1,\ldots,q_n$. Then $(w_1,\ldots,w_n)$ is
periodic with period~$\frac{2\pi r}{A\mu}$.

In case (ii), since $u$ is periodic with period $S$, if
$(w_1,\ldots,w_n)$ is periodic with period $T$ then $T=rS$ for
some $r\in\N$. But then $e^{i\phi_j(s+T)}=e^{i\phi_j(s)}$ for
$j=1,\ldots,n$, so \eq{lm6eq3} gives $e^{ir\ga_j}=1$, and
$\ga_j\in 2\pi\Z/r\subset\pi\Q$ for $j=1,\ldots,n$, as we want.
Conversely, if $\ga_j\in\pi\Q$ for $j=1,\ldots,n$ then we may
write $\ga_j=2\pi p_j/r$ for $j=1,\ldots,n$, $p_j\in\Z$ and
$r\in\N$, and then $(w_1,\ldots,w_n)$ is periodic with
period~$T=rS$.

It remains to show that in both cases, for fixed $m,\al$, there is
a dense subset of initial data with $(w_1,\ldots,w_n)$ periodic.
In case (i) this is straightforward: by Theorem E, for fixed
$\la_1,\ldots,\la_n$ and $\al$, solutions are in 1-1
correspondence with choices of $\al_1,\ldots,\al_n>0$ and
$\psi_1,\ldots,\psi_n \in\R$ satisfying
$\sum_{j=1}^n\frac{\la_j}{\al_j}+\al=0$, and by the previous part,
the corresponding solution is periodic if and only if
$\smash{\frac{\la_j}{\al_j}}=\mu q_j$ for $\mu>0$ and $q_j\in\Q$
for $j=1,\ldots,n$. It is easy to see that the set of such
$\al_j,\psi_j$ is dense in the set of all allowed~$\al_j,\psi_j$.

So we restrict to case (ii). In the special Lagrangian case
$\al=0$, the first author \cite[\S 5.5]{Joyc1} showed that
periodic solutions are dense in all solutions, so we suppose
$\al\ne 0$. Given some solution in Theorem E, Lemma \ref{lm6lem1}
found a unique $u^*\in[u_1,u_2]\subset (\be_1,\be_2)$ where $G(u)$
is maximum in $(\be_1,\be_2)$, and Theorem E showed that
$u:\R\ra\R$ cycles between $u_1$ and $u_2$ and so realizes all
values in $[u_1,u_2]$, including $u^*$. Thus, in Theorem B we can
choose the base point $s_0\in I=\R$ so that $u(s_0)=u^*$;
effectively, this changes $\al_j\mapsto \al_j+\la_ju^*$, $u\mapsto
u-u^*$, and~$u^*\mapsto 0$.

We will work for the rest of the proof with this normalization, so
that $u^*=0$. Then $G(u)$ has a maximum at $u=0$, so that
\eq{lm5eq5} gives
\begin{equation}
\ts\sum_{j=1}^n\frac{\la_j}{\al_j}+\al=0.
\label{lm6eq5}
\end{equation}

The remaining variables are $\al_1,\ldots,\al_n>0$ which satisfy
\eq{lm6eq5}, $A$ which satisfies $0<A<(\al_1\cdots\al_n)^{1/2}$ by
\eq{lm4eq9} and (ii), and $\psi_1,\ldots,\psi_n\in\R$. Now
$\ga_1,\ldots,\ga_n$ are independent of $\psi_1,\ldots,\psi_n$, so
we can regard them as functions of $\al_1,\ldots,\al_n$ and $A$.
Define
\begin{align}
&\Psi^{m,n}:\bigl\{(\al_1,\ldots,\al_n,A)
\in(0,\iy)^{n+1}:
\nonumber\\
&\ts\sum_{j=1}^n\frac{\la_j}{\al_j}+\al=0,\;\>
A<(\al_1\cdots\al_n)^{1/2}\bigr\}\longra\R^n,
\quad\text{where}
\label{lm6eq6}\\
&\Psi^{m,n}=(\Psi^{m,n}_1,\ldots,\Psi^{m,n}_n):(\al_1,\ldots,\al_n,A)
\longmapsto (\ga_1,\ldots,\ga_n). \nonumber
\end{align}

To compute $\Psi^{m,n}_j$ explicitly, note that in one period $S$ of
$s$, $u$ goes from $u_1$ up to $u_2$ and back down again, and
$\psi_j$ increases by $\frac{\ga_j}{2}$ in each half-period. So
changing variables from $s$ to $u$ in $[u_1,u_2]$ we see that
$\ga_j=2\int_{u_1}^{u_2}\frac{\d\phi_j}{\d u}(u)\d u$, taking the
branch of $\frac{\d\phi_j}{\d u}$ on $(u_1,u_2)$ for which $\frac{\d
u}{\d s}>0$. Computing $\frac{\d\phi_j}{\d u}$ from \eq{lm4eq7} and
using \eq{lm4eq9} to eliminate terms in
$\sin(\phi-\th),\cos(\phi-\th)$ gives
\begin{equation}
\Psi^{m,n}_j(\al_1,\ldots,\al_n,A)=-\int_{u_1}^{u_2}\frac{A\la_j\,\d
v}{(\al_j+\la_jv)\sqrt{Q(v)e^{\al v}-A^2}}\,, \label{lm6eq7}
\end{equation}
where $u_1<u^*=0<u_2$ are the closest roots of $Q(v)e^{\al v}=A^2$
to zero.

We must prove that $\Psi^{m,n}(\al_1,\ldots,\al_n,A)\in(\pi\Q)^n$
for a dense subset of $(\al_1,\ldots,\al_n,A)$ in the domain of
$\Psi^{m,n}$. To do this we will use the method of Joyce \cite[\S
5.5]{Joyc1}. We first compute various limits of~$\Psi^{m,n}$.

\begin{prop} Regarding $\al_1,\ldots,\al_n>0$ satisfying \eq{lm6eq5}
as fixed, for all\/ $j$ we have
\begin{equation}
\lim_{A\ra(\al_1\cdots\al_n)^{1/2}_-\!\!} \Psi^{m,n}_j
(\al_1,\ldots,\al_n,A)\!=\!\ts-2\pi\la_j\al_j^{-1}\bigl(
2\sum\limits_{k=1}^n\la_k^2 \al_k^{-2}\bigr)^{-1/2}.
\label{lm6eq8}
\end{equation}
\label{lm6prop2}
\end{prop}

\begin{proof} Recall that $(\al_1\cdots\al_n)^{1/2}\sin(\phi-\th)=
A$ at $u=0$ from \eq{lm3eq6}. When $A$ is close to
$(\al_1\cdots\al_n)^{1/2}$, $u$ is small and $\sin(\phi-\th)$ is
close to 1, so $\phi-\th$ remains close to $\pi/2$. Write
$\phi-\th=\frac{\pi}{2}+\varphi$, for $\varphi$ small. Then, setting
$Q(u)\approx\al_1\cdots\al_n$,
\begin{equation*}
\cos(\phi-\th)\approx -\varphi,\;\> \sin(\phi-\th)\approx 1
\;\>\text{and}\;\> \sum_{k=1}^n\frac{\la_k}{\al_k+\la_ku}+\al
\approx -u\sum_{k=1}^n\la_k^2\al_k^{-2}
\end{equation*}
via linear approximation, taking only the highest order terms,
equation \eq{lm4eq7} implies that
\begin{align*}
& \frac{\d u}{\d s}\approx -2(\al_1\cdots\al_n)^{1/2}\varphi,\quad
\frac{\d(\phi-\th)}{\d s}=\frac{\d\varphi}{\d s}\approx
u(\al_1\cdots\al_n)^{1/2}\sum_{k=1}^n\la_k^2\al_k^{-2}.
\end{align*}
It follows that
\begin{align*} & \frac{\d^2 u}{\d s^2}+ \Bigl(2\al_1\cdots\al_n
\sum_{k=1}^n\la_k^2\al_k^{-2}\Bigr)u \approx 0 \qquad \text{and}\\
&\frac{\d^2 (\phi-\th)}{\d s^2}+ \Bigl(2\al_1\cdots\al_n
\sum_{k=1}^n\la_k^2\al_k^{-2}\Bigr)(\phi-\th) \approx 0,
\end{align*}
so that $u$ and $\phi-\th$ undergo approximately simple harmonic
oscillations with period
$S=2\pi\bigl(2\al_1\cdots\al_n\sum_{k=1}^n\la_k^2\al_k^{-2}
\bigr)^{-1/2}$. Then \eq{lm4eq7} shows that
\begin{equation*}
\frac{\d\phi_j}{\d s}\approx
-\la_j\al_j^{-1}(\al_1\cdots\al_n)^{1/2},
\end{equation*}
which is approximately constant. Hence
\begin{equation*}
\ga_j=\int_0^S\frac{\d\phi_j}{\d s}\d s\approx\frac{\d\phi_j}{\d
s}S=-2\pi\la_j\al_j^{-1}\Bigl( 2\sum_{k=1}^n\la_k^2
\al_k^{-2}\Bigr)^{-1/2}
\end{equation*}
 This proves~\eq{lm6eq8}.
\end{proof}

For an inductive step needed later, we have to compute what happens
when $\al_n\ra\iy$ or $\al_1\ra\iy$ and include the case $m=0$ when
$\al>0$, which is also well-defined using \eq{lm6eq6} and
\eq{lm6eq7}. So we allow $1\le m\le n$ when $\al\leq 0$, and $0\le
m<n$ when~$\al>0$.

\begin{prop} Suppose $\bigl(\al_1(t),\ldots,\al_n(t),A(t)\bigr),$
$t\in(1,\iy),$ is a continuous path in the domain of\/ $\Psi^{m,n}$
in {\rm\eq{lm6eq6},} such that
\begin{equation}
\begin{gathered} \text{$\lim_{t\ra\iy}\al_j(t)=\ti\al_j$\quad for
$j=1,\ldots,n-1,$
\qquad $\lim_{t\ra\iy}\al_n(t)=\iy,$} \\
\text{and\/\qquad $\lim_{t\ra\iy}A(t)\al_n(t)^{-1/2}=\ti A.$}
\end{gathered}
\label{lm6eq9}
\end{equation}
Then $(\ti\al_1,\ldots,\ti\al_{n-1},\ti A)$ is in the domain of\/
$\Psi^{\ti m,n-1},$ with\/ $\ti m=\min(m,\ab n-1),$ and
\begin{equation*}
\lim_{t\ra\iy}\!\Psi^{m,n}_j\bigl(\al_1(t),\ldots,\al_n(t),A(t)\bigr)
\!=\!\begin{cases}
\Psi^{\ti m,n-1}_j(\ti\al_1,\ldots,\ti\al_{n-1},\ti A), & j\!<\!n, \\
0, & j\!=\!n.\end{cases}
\end{equation*}

Similarly, suppose $\bigl(\al_1(t),\ldots,\al_n(t),A(t)\bigr),$
$t\in(1,\iy),$ is a continuous path in the domain of\/ $\Psi^{m,n}$
such that
\begin{gather*}
\text{$\lim_{t\ra\iy}\al_j(t)=\ti\al_j$\quad for $j=2,\ldots,n,$
\qquad $\lim_{t\ra\iy}\al_1(t)=\iy,$} \\
\text{and\/\qquad $\lim_{t\ra\iy}A(t)\al_1(t)^{-1/2}=\ti A.$}
\end{gather*}
Then $(\ti\al_2,\ldots,\ti\al_n,\ti A)$ is in the domain of\/
$\Psi^{\ti m-1,n-1},$ with\/ $\ti m\!=\!\max(m,1),$ and
\begin{equation*}
\lim_{t\ra\iy}\!\Psi^{m,n}_j\bigl(\al_1(t),\ldots,\al_n(t),A(t)\bigr)
\!=\!\begin{cases} \Psi^{\ti m-1,n-1}_{j-1}(\ti\al_2,\ldots,
\ti\al_n, \ti A), & j\!>\!1, \\ 0, & j\!=\!1.\end{cases}
\end{equation*}
\label{lm6prop3}
\end{prop}

\begin{proof} We have $Q(v)=\prod_{j=1}^n(\al_j(t)+\la_jv)$. Write
$\ti Q(v)=\prod_{j=1}^{n-1}(\ti\al_j+\la_jv)$. Then in the case of
\eq{lm6eq9}, we see that the integrand in \eq{lm6eq7} satisfies
\begin{align*}
&\lim_{t\ra\iy}\frac{A\la_j}{(\al_j+\la_jv)\sqrt{Q(v)e^{\al
v}-A^2}}\\
&=\!\lim_{t\ra\iy}\frac{1}{\al_j(t)\!+\!\la_jv} \lim_{t\ra\iy}
\frac{A(t)(\al_n(t)+\la_nv)^{-1/2}\la_j}{
\sqrt{\prod_{l=1}^{n-1}(\al_l(t)\!+\!\la_lv)e^{\al
v}\!-\!A(t)^2(\al_n(t)\!+\!\la_nv)^{-1}}}\\
&=\begin{cases}\frac{\ti
A\la_j}{\sqrt{\prod_{l=1}^{n-1}(\ti\al_l+\la_lv)e^{\al v}-\ti
A^2}}\cdot\frac{1}{\ti\al_j+\la_jv}, & j=1,\ldots,n-1,\\
0, & j=n.\end{cases}
\end{align*}
Applying the Dominated Convergence Theorem to \eq{lm6eq7}, and
noting that $u_1(t)\ra\ti u_1$, $u_2(t)\ra\ti u_2$, gives
\begin{align*}
\lim_{t\ra\iy}\Psi^{m,n}_j\bigl(\al_1(t),\ldots,\al_n(t),A(t)\bigr)&=
-\int_{\ti u_1}^{\ti u_2}\frac{\ti A\la_j\,\d
v}{(\ti\al_j+\la_jv)\sqrt{\ti Q(v)e^{\al v}-\ti A^2}}\\
&=\Psi^{\ti m,n-1}_j(\ti\al_1,\ldots,\ti\al_{n-1},\ti A)
\end{align*}
if $j=1,\ldots,n-1$, and shows that the limit is 0 if $j=n$. This
proves the first part of the proposition, and the second part is
similar.
\end{proof}

Note that when $\al>0$ and $m=n-1$, from \eq{lm6eq5} we always have
$\al_n$ bounded. Similarly, when $\al<0$ and $m=1$, we always have
$\al_1$ bounded. Hence the first and second parts of Proposition
\ref{lm6prop3}, respectively, cannot apply. In the following, we
only need the first part of Proposition \ref{lm6prop3} when $\al<0$,
and the second part when~$\al>0$.

\begin{prop} For $\Psi^{m,n}$ as in \eq{lm6eq6} and\/ {\rm\eq{lm6eq7},}
where we allow $m=n$ only if\/ $\al<0$ and\/ $m=0$ only if\/
$\al>0,$ the image $\Image\Psi^{m,n}$ is $n$-dimensional, and for a
dense open subset of\/ $(\al_1,\ldots,\al_n,A)$ in the domain of\/
$\Psi^{m,n},$ the following derivative is an isomorphism:
\begin{equation}
\d\Psi^{m,n}\vert_{(\al_1,\ldots,\al_n,A)}:\bigl\{(x_1,\ldots,x_n,y)
\!\in\!\R^{n+1}:\ts\sum_{j=1}^n\frac{\la_jx_j}{\al_j^2}\!=\!0\bigr\}
\!\ra\!\R^n,
\label{lm6eq13}
\end{equation}
so that\/ $\Psi^{m,n}$ is a local diffeomorphism
near~$(\al_1,\ldots,\al_n,A)$.
\label{lm6prop4}
\end{prop}

\begin{proof} By Proposition \ref{lm6prop2} the closure
$\ov{\Image\Psi^{m,n}}$ contains the set
\begin{equation*}
\left\{\left(\frac{-2\pi\la_1\al_1^{-1}}{( 2\sum_{k=1}^n\la_k^2
\al_k^{-2})^{1/2}},\cdots, \frac{-2\pi\la_n\al_n^{-1}}{(
2\sum_{k=1}^n\la_k^2 \al_k^{-2})^{1/2}}\right):\text{$\al_ j>0$ for
all $j$}\right\},
\end{equation*}
which is a nonempty open subset of the $(n-1)$-dimensional real
hypersurface
\begin{equation*}
H=\bigl\{(\ga_1,\ldots,\ga_n)\in\R^n:\ts\sum_{j=1}^n\ga_j^2=
2\pi^2\bigr\}
\end{equation*}
in $\R^n$. This implies that $\Image\Psi^{m,n}$ is at least
$(n-1)$-dimensional. Since $\Psi^{m,n}$ is real analytic and its
domain is nonsingular and connected, there are only two
possibilities:
\begin{itemize}
\setlength{\itemsep}{0pt}
\setlength{\parsep}{0pt}
\item[(A)] $\Image\Psi^{m,n}$ is $n$-dimensional, or
\item[(B)] $\Image\Psi^{m,n}$ lies in the $(n-1)$-dimensional real
hypersurface $H$ in~$\R^n$.
\end{itemize}

We shall use Proposition \ref{lm6prop3} and induction on $n$ to
eliminate possibility (B), so that (A) holds. The first step $n=1$
is studied by Abresch and Langer \cite{AbLa}, and translated into
our notation, \cite[Th.~A \& Prop.~3.2(v)]{AbLa} implies that when
$\al<0$ and $m=n=1$, $\Image\Psi^{1,1}=(-\sqrt{2}\pi,-\pi)$.
Changing signs of $\la_1,\al$ we deduce that when $\al>0$, $m=0$ and
$n=1$, $\Image\Psi^{0,1}=(\pi,\sqrt{2}\pi)$. Thus in both cases
$\Image\Psi^{m,n}$ is $n$-dimensional when~$n=1$.

Suppose by induction that $n\ge 2$, and that $\Image\Psi^{k,l}$ is
$l$-dimensional whenever $l<n$, allowing $k=l$ only if $\al<0$, and
$k=0$ only if $\al>0$. First suppose $\al<0$, and let $0<m\le n$.
Set $\ti m=\min(m,n-1)$. Then all of the domain of $\Psi^{\ti
m,n-1}$ arises as limits of the domain of $\Psi^{m,n}$ as in
\eq{lm6eq9}, so the first part of Proposition \ref{lm6prop3} implies
that $\Image\Psi^{\ti m,n-1}\t\{0\}\subset\ov{\Image\Psi^{m,n}}$.
But by induction $\Image\Psi^{\ti m,n-1}$ is $(n-1)$-dimensional, so
$\Image\Psi^{\ti m,n-1}\t\{0\}$ is not contained in the hypersurface
$H$ in $\R^n$ because $\bigl(\Image\Psi^{\ti m,n-1}\t\{0\}\bigr)
\cap H$ is at most $(n-2)$ dimensional. So (B) does not hold.

Similarly, if $\al>0$ then for $0\le m<n$ and $\ti m=\max(m,1)$, the
second part of Proposition \ref{lm6prop3} implies that
$\{0\}\t\Image\Psi^{\ti m-1,n-1}\subset\ov{\Image\Psi^{m,n}}$, and
$\{0\}\t\Image\Psi^{\ti m-1,n-1}$ is $(n-1)$-dimensional and not
contained in $H$, so (B) does not hold. Thus in both cases (A)
holds, so $\Image\Psi^{m,n}$ is $n$-dimensional, proving the
inductive step. The final parts follow as $\Psi^{m,n}$ is real
analytic and its domain is nonsingular and connected.
\end{proof}

We can now complete the proof of Theorem F. By Proposition
\ref{lm6prop4}, $\Psi^{m,n}$ is a local diffeomorphism near
$(\al_1,\ldots,\al_n,A)$ for $(\al_1,\ldots,\al_n,A)$ in a dense
open subset $U$ in the domain of $\Psi^{m,n}$. As $(\pi\Q)^n$ is
dense in the range $\R^n$ of $\Psi^{m,n}$, it follows that
$(\Psi^{m,n})^{-1}\bigl((\pi\Q)^n\bigr)$ is dense in $U$, and hence
in the domain of $\Psi^{m,n}$, since $U$ is dense. But a choice of
initial data gives a periodic solution if and only if
$\Psi^{m,n}(\al_1,\ldots,\al_n,A)\in(\pi\Q)^n$, by a previous part
of the theorem. Since this holds for a dense subset of allowed
$(\al_1,\ldots,\al_n,A)$, a dense subset of choices of initial data
yield periodic solutions.
\end{proof}

\end{document}